# ON THE ŁOJASIEWICZ EXPONENT OF THE GRADIENT OF A POLYNOMIAL FUNCTION

Andrzej Lenarcik

*Let nothing be done through selfish ambition or conceit (Flp 2,3)*

## Abstract

Let $h = \sum h_{\alpha\beta} X^\alpha Y^\beta$ be a polynomial with complex coefficients. The Łojasiewicz exponent of the gradient of $h$ at infinity is the upper bound of the set of all real $\lambda$ such that $|\operatorname{grad} h(x,y)| \geq c|(x,y)|^\lambda$ in a neighbourhood of infinity in $\mathbf{C}^2$, for $c > 0$. We estimate this quantity in terms of the Newton diagram of $h$. The equality is obtained in the nondegenerate case.

## 1  Introduction

The Łojasiewicz exponent $l_\infty(H)$ of the polynomial mapping $H = (f, g) : \mathbf{C}^2 \mapsto \mathbf{C}^2$ at infinity is the upper bound of the set of all real $\lambda$ such that

$$|H(z)| \geq c|z|^\lambda \qquad (1)$$

for sufficiently large $|z|$ and for $c > 0$. If the set of all the exponents is empty we put $l_\infty(H) = -\infty$. We use the norm $|z| = \max\{|x|, |y|\}$ for $z = (x, y) \in \mathbf{C}^2$. The quantity $l_\infty(H)$ is also called the exponent of growth of a polynomial mapping $H$.

Chądzyński and Krasiński [ChK1] showed that the number of solutions of the equation $f = g = 0$ is finite if, and only if, $l_\infty(f, g) > -\infty$, and then the exponent is realized on at least one of the curves $\{f = 0\}$ or $\{g = 0\}$. They also proved that $l_\infty(H)$ is a rational number or $-\infty$. In the paper [ChK2] the authors described $l_\infty(H)$ using resultant. In [Pł1] Płoski gave an estimation of $l_\infty(H)$ for a polynomial mapping $\mathbf{C}^n \mapsto \mathbf{C}^n$ in terms of its geometrical degree and the degrees of the mapping components. The properness of $H$ can be characterized by using $l_\infty(H)$ ($H$ is proper iff $l_\infty(H) > 0$). This exponent is also applicable in the theory of polynomial authomorphisms, especially the exponent of gradient. Chądzyński and Krasiński showed in [ChK2] that a polynomial $h : \mathbf{C}^2 \to \mathbf{C}$ is the component of a polynomial automorphism if, and only if, the system of equations $\frac{\partial h}{\partial X} = \frac{\partial h}{\partial Y} = 0$ has no solutions and $l_\infty(\operatorname{grad} h) > -1$ [ChK2]. A connection of $l_\infty(\operatorname{grad} h)$ with the Newton diagram for nondegenerate $h$ were observed by Pi. Cassou-Noguès and Há Huy Vui [CN-H]. It seems to be some imprecisions in the formulation of Proposition 10 (page 42). For $h(X,Y) = Y^p + X^p$ (page 24) $l_\infty(\operatorname{grad} h) = p - 1$, but by the proposition this exponent equals $p$. The aim of our paper is to give an estimation of $l_\infty(\operatorname{grad} h)$ in terms of the Newton diagram of $h$, without the nondegeneracy assumption. For nondegenerate polynomials, the equality is obtained. Our methods are different from the methods used in [CN-H]. The results presented in our paper are the counterparts of the results obtained by the author in the local case [L].

We consequently use conventions: $\inf \emptyset = +\infty$ and $\sup \emptyset = -\infty$.

## 2 The main result

A presentation of the main result needs some definitions. Let $h(X,Y) = \sum h_{\alpha\beta} X^\alpha Y^\beta$ be a polynomial with complex coefficients. We define the *support* of $h$ as the set $\{(\alpha, \beta) : h_{\alpha\beta} \neq 0\}$ and denote it by $\operatorname{supp} h$. The *degrees* $\deg h$, $\deg_X h$ and $\deg_Y h$ are defined to be the maximae of the expressions $\alpha + \beta$, $\alpha$ and $\beta$, respectively, where $(\alpha, \beta)$ runs over $\operatorname{supp} h$. For $h = 0$ we put $-\infty$ for each of the above degrees. Analogously, we define the *orders* $\operatorname{ord} h$, $\operatorname{ord}_X h$ and $\operatorname{ord}_Y h$ as the minimae of the respective expressions. For $h = 0$ we put $+\infty$ for each order. The *Newton diagram* $\Delta_h$ of $h$ is the convex hull of $\operatorname{supp} h$. The *set of boundary segments* is the set of all one-dimensional faces which compose the boundary of $\Delta_h$. For any boundary segment $S$ we define $\operatorname{in}(h, S)$ as the sum of monomials $h_{\alpha\beta} x^\alpha y^\beta$ over all $(\alpha, \beta) \in S$. We say that $h$ is *nondegenerate* on $S$ if the system of equations $\frac{\partial}{\partial X}\operatorname{in}(h, S) = \frac{\partial}{\partial Y}\operatorname{in}(h, S) = 0$ has no solutions in $(\mathbf{C} \setminus \{0\}) \times (\mathbf{C} \setminus \{0\})$. The *right Newton polygon* $\mathcal{N}_h^{(\mathbf{r})}$ consists of all the boundary segments which lay on the right side of $\Delta_h$ and join the lines $\beta = \operatorname{ord}_Y h$ and $\beta = \deg_Y h$. Analogously, the *top Newton polygon* $\mathcal{N}_h^{(\mathbf{t})}$ consists of all the segments which lay on the top of $\Delta_h$ and join the lines $\alpha = \operatorname{ord}_X h$ and $\alpha = \deg_X h$. The set of segments of the both polygons is called the *Newton polygon of $h$ at infinity*. If $h$ is nondegenrate on each segment of this polygon, then we say that $h$ is *nondegenerate at infinity*.

Usually, the Newton Polygon at infinity of a polynomial $h$ is defined to be the set of all the boundary segments, not included in the axes, for $h$+generic const. This definition coincides with the earliest one if $h(X,0)h(0,Y) \neq 0$, but differs in general (e.g. $h = X + X^2 Y$).

We say that a segment of $\mathcal{N}_h^{(\mathbf{r})}$ is *exceptional* if it joins the horizontal axis with a point of the form $(p, 1)$. Analogously, a segment of $\mathcal{N}_h^{(\mathbf{t})}$ is exceptional if it joins the vertical axis with a point of the form $(1, q)$. Notice, that the right polygon has no more than one exceptional segment.

Similarly — top polygon. It is convenient to put the segments of each of the polygons $\mathcal{N}_h^{(\mathbf{r})}, \mathcal{N}_h^{(\mathbf{t})}$ in order. We define an order in the right polygon in such a way that the first segment is the nearest to the horizontal axis. An analogous order in the top polygon is defined such that the first segment is the nearest to the vertical axis. Notice, that if the exceptional segment of the right or top polygon exists, then it is the first segment of the polygon.

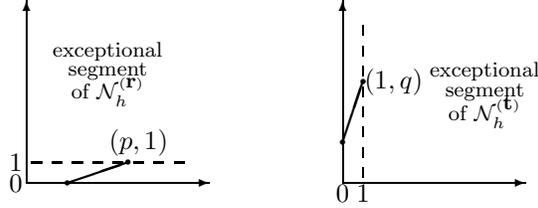

EXAMPLE. Let $h = X^2 + X^4 + XY^3 + XY^6 + X^7Y + X^4Y^8 + X^9Y^4 + X^9Y^6 + X^7Y^8$. The Newton polygon $\Delta_h$ has the nine boundary segments $A$, $B$, $C$, $D$, $E$, $F$, $G$, $H$, $I$ which join the vertices $(2,0)$, $(4,0)$, $(7,1)$, $(9,4)$, $(9,6)$, $(7,8)$, $(4,8)$, $(1,6)$, $(1,3)$, $(2,0)$, respectively. We have $\mathcal{N}_h^{(\mathbf{r})} = \{B, C, D, E\}$ and $\mathcal{N}_h^{(\mathbf{t})} = \{G, F, E\}$. The polygon of $h$ at infinity is composed of all the boundary segments except $A$, $H$ and $I$. Notice, that $B$ is the first segment of the right polygon, and $G$ is the first segment of the top polygon, however, $B$ is exceptional in $\mathcal{N}_h^{(\mathbf{r})}$.

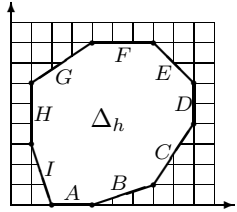

For a segment $S$, not parallel to the horizontal axis, we denote by $\alpha(S)$ the abscissa of the point, where the line determined by $S$ intersects the horizontal axis. Analogously, for $S$, not parallel to the vertical axis, we denote by $\beta(S)$ the ordinate of the point, where the line intersects the vertical axis. The following theorem is the main result of the paper.

THEOREM 2.1. *Let $h \in \mathbf{C}[X, Y]$ be a polynomial without constant term, with non-zero gradient's components, not divisible by $X^2$ and not divisible by $Y^2$. If, additionally, at least one of the polygons $\mathcal{N}_h^{(\mathbf{r})}$, $\mathcal{N}_h^{(\mathbf{t})}$ has a segment which is not exceptional, then*

$$l_\infty(\operatorname{grad} h) \leq \min \left\{ \inf_{S'} \alpha(S'), \inf_{S''} \beta(S'') \right\} - 1 ,$$

*where $S'$ runs over all the segments of the right polygon without exceptional one, and $S''$ runs over all the segments of the top polygon without exceptional one. Moreover, if $h$ is nondegenerate at infinity, then the equality holds.*

A proof of the theorem is given in Section 9. Each of the infimae in the statement means minimum, if the corresponding set of segments is nonempty, and $+\infty$, otherwise. Notice, that the assumption $h(0,0) = 0$ is not restrictive. We can always consider $h - h(0,0)$ with the same gradient. If one of the gradient's components vanish then $l_\infty(\operatorname{grad} h) = 0$ or $-\infty$. If $h$ is divisible by $X^2$ or $Y^2$, then, obviously, $l_\infty(\operatorname{grad} h) = -\infty$. The case, when the both considered sets of segments are empty, is described in the following elementary

PROPOSITION 2.2. *Let $h \in \mathbf{C}[X,Y]$ be a polynomial without constant term, with non-zero gradient's components, not divisible by $X^2$ and not divisible by $Y^2$. If each of the polygons $\mathcal{N}_h^{(\mathbf{r})}$ and $\mathcal{N}_h^{(\mathbf{t})}$ is empty or has an exceptional segment only, then*

$$h(X,Y) = aX + bY + cXY ,$$

*where $ab \neq 0$ or $c \neq 0$. (For a proof, see Propositions 9.1 and 9.2).*

Under assumptions of the proposition we have $l_\infty(\operatorname{grad} h) = 1$ if $c \neq 0$ and $l_\infty(\operatorname{grad} h) = 0$ if $c = 0$, $ab \neq 0$.

A nondegeneracy of a polynomial $h$ on a boundary segment $S$ can be easily examined. If the line determined by $S$ does not intersect the origin, then the nondegeneracy of $h$ on $S$ is equivalent to the fact that $\operatorname{in}(h,S)$ has no multiple factors different from $X$ and $Y$. If the line determined by $S$ intersects the origin, and $(0,0)$ is not the end of the segment, then $h$ is degenerate on $S$.

Let us return to the polynomial $h$ considered in the example before Theorem 2.1. It is easy to verify that all the assumptions of the theorem are satisfied. Moreover, $h$ is nondegenerate. The infimae from the statement of the theorem are attained for the earliest segments, in the sense of order in the polygons. The first segment of right polygon, which is not exceptional, is $C$. Similarly, for the top polygon, it is $G$. Hence, $l_\infty(\operatorname{grad} h) = \min\{\alpha(C), \beta(G)\} - 1 = \min\{6\frac{1}{3}, 5\frac{1}{3}\} - 1 = 4\frac{1}{3}$. The necessity of omission the exceptional segments is not observed in [CN-H].

As in the local case [L], the proof of Theorem 2.1 is based on the formula which describes the Łojasiewicz exponent of a pair of polynomials by using information from the Newton diagrams of the both components. This is the subject of the next section.

## 3 Auxiliary results

Let $H = (f,g)$ be a pair of polynomials. Our aim, in this section, is to describe a connections between $l_\infty(H)$ and the Newton diagrams $\Delta_f$ and $\Delta_g$ (Theorems 3.1 and 3.2). This results are used in the proof of Theorem 2.1 in Section 9.

First, we give some definitions. For a segment $S$ we denote by $S_1$ and $S_2$ the projections of $S$ on the horizontal and vertical axis, respectively, and by $|S_1|$ and $|S_2|$ their lengths. We define a number $\sigma(S) \in \{-1, 0, 1\}$. If $S$ is parallel to one of the axes we put $\sigma(S) = 0$. In other cases $\sigma(S)$ has an opposite value to the sign of the slope of $S$. In the example before Theorem 2.1 we have $\sigma = 0$ for $A$, $D$, $F$, $H$, $\sigma = 1$ for $E$, $I$, and $\sigma = -1$ for $B$, $C$, $G$. Consider a polynomial and its global Newton polygon. The *declivity of a segment $S$ of the right Newton polygon* is defined to be the number $\frac{|S_1|}{|S_2|}\sigma(S)$. Obviously, it is well defined, and considering the order of the segments in the polygon, it is an increasing function. Analogously, a *declivity of a segment $S$ of the top Newton polygon* is defined to be the number $\frac{|S_2|}{|S_1|}\sigma(S)$. It is also well defined and increasing in the same sense. In the sequel, we write simply "the declivity of $S$" if it is clear what polygon (right or top) the segment belongs to. For any non-zero $h$ and for any segment $S$ such that $|S_2| \neq 0$ we define a number

$$\alpha(S, \Delta_h) = \max\left\{\alpha + \beta\frac{|S_1|}{|S_2|}\sigma(S) : (\alpha, \beta) \in \operatorname{supp} h\right\}, \qquad (2)$$

and for $S$ such that $|S_1| \neq 0$ a number

$$\beta(S, \Delta_h) = \max\left\{\alpha\frac{|S_2|}{|S_1|}\sigma(S) + \beta : (\alpha, \beta) \in \operatorname{supp} h\right\}. \qquad (3)$$

These numbers have simple geometrical meaning. The first is equal to the maximal possible abscissa of the point where the line supporting $\Delta_h$, parallel to $S$, intersects the horizontal axis. The second number is equal to the maximal possible ordinate of the point where the line of the same type intersects the vertical axis.

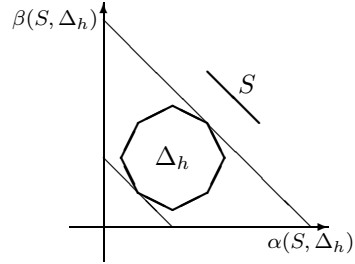

The above-introduced numbers allow as to give an upper estimation of $l_\infty(H)$ by using Newton diagrams $\Delta_f$ and $\Delta_g$. We have the following

THEOREM 3.1. *For a pair $H = (f, g)$ of non-zero polynomials the exponent $l_\infty(H)$ is bounded from above by the minimum of the following six quantities:*

$$\deg H(X, 0), \quad \inf_{S \in \mathcal{N}_f^{(\mathbf{r})}} \alpha(S, \Delta_g), \quad \inf_{T \in \mathcal{N}_g^{(\mathbf{r})}} \alpha(T, \Delta_f),$$

$$\deg H(0,Y), \quad \inf_{S \in \mathcal{N}_f^{(\mathbf{t})}} \beta(S, \Delta_g), \quad \inf_{T \in \mathcal{N}_g^{(\mathbf{t})}} \beta(T, \Delta_f),$$

*where, by the degree of a pair we understand the maximum of its components' degrees.*

As earlier, we use infimae to preserve a sense of the quantities when the corresponding polygons are empty. The proof of the theorem is given in Section 7. Analogously, as in Theorem 2.1 we may obtain equality in the nondegenerate case. First, we must give a suitable

DEFINITION. *We say that a pair of non-zero polynomials $H = (f,g)$ is nondegenerate at infinity if for any segment $S$ of the polygon of $f$ at infinity and for any segment $T$ of the polygon of $g$ at infinity one of the following conditions holds:*

(a) *$S$ and $T$ are not parallel,*
(b) *$S$ and $T$ are parallel, $S \in \mathcal{N}_f^{(\mathbf{r})}$, $T \in \mathcal{N}_g^{(\mathbf{r})}$ and the system of equations $in(f,S) = in(g,T) = 0$ has no solutions in $(\mathbf{C} \setminus \{0\}) \times (\mathbf{C} \setminus \{0\})$.*
(c) *$S$ and $T$ are parallel, $S \in \mathcal{N}_f^{(\mathbf{t})}$, $T \in \mathcal{N}_g^{(\mathbf{t})}$ and the system of equations $in(f,S) = in(g,T) = 0$ has no solutions in $(\mathbf{C} \setminus \{0\}) \times (\mathbf{C} \setminus \{0\})$.*

We have the following

THEOREM 3.2. *If $H = (f,g)$ is a pair of non-zero polynomials, nondegenerate at infinity, then the Lojasiewicz exponent $l_\infty(H)$ is equal to the minimum of the six quantities given in Theorem 3.1.*

The proof of the theorem is given in Section 7. Theorems 3.1 and 3.2 have its local counterpart ([L], Theorem 4.2).

EXAMPLE. Let $H = (f,g)$, where $f = Y^2 + X^4Y^4 + X^5Y^7 + X^3Y^8$ and $g = X^2 + X^3 + X^7Y + X^6Y^4$. The polygon of $f$ at infinity consists of the four segments $A$, $B$, $C$, $D$ which join the vertices $(0,2)$, $(4,4)$, $(5,7)$, $(3,8)$, $(0,2)$, respectively. We have $\mathcal{N}_f^{(\mathbf{r})} = \{A,B,C\}$ and $\mathcal{N}_f^{(\mathbf{t})} = \{D,C\}$. The polygon of $g$ at infinity consists of the three segments $E$, $F$, $G$ which join the vertices $(3,0)$, $(7,1)$, $(6,4)$, $(2,0)$, respectively. We have $\mathcal{N}_g^{(\mathbf{r})} = \{E,F\}$ and $\mathcal{N}_g^{(\mathbf{t})} = \{G,F\}$.

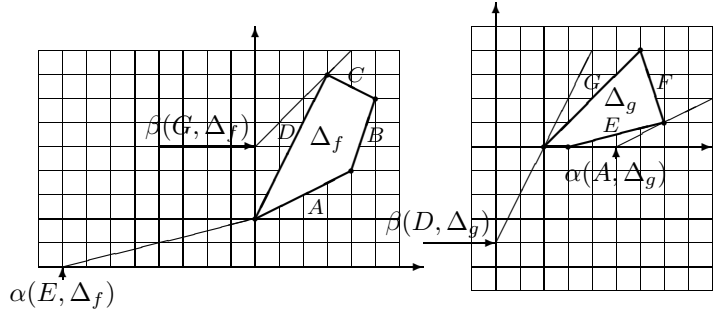

All the considered polygons are non-empty. In this case, the first two infimae, in the statement of Theorem 3.1, are attained for the segments nearest to the horizontal axis. The next two infimae are attained for the segments nearest to the vertical axis. So, the six quantities from Theorem 3.1 are respectively equal to: $\deg H(X,0) = 3$, $\alpha(A, \Delta_g) = 5$, $\alpha(E, \Delta_f) = -8$, $\deg H(0,Y) = 2$, $\beta(D, \Delta_g) = -4$, $\alpha(G, \Delta_f) = 5$. Obviously, the pair is non-degenerate. From Theorem 3.2 we have $l_\infty(H) = \min\{3, 5, -8, 2, -4, 5\} = -8$.

## 4 Relative exponents

In this short section we define a few versions of the Łojasiewicz exponent related to a variable or to a subset. They simplify the process of calculating the Łojasiewicz exponent of a polynomial mapping at infinity.

Let $H = (f, g)$ be a pair of polynomials. Analogously, as in the local case ([L], [Pł2]) we can consider a relative exponent $l_\infty(H, X)$ which is defined to be the upper bound of the set of all real $\lambda$ such that inequality

$$|H(x,y)| \geq c|x|^\lambda \tag{4}$$

holds for sufficiently large $|x|$ and for $c > 0$. The exponent $l_\infty(H, Y)$ is defined analogously. In the local situation, it is easy to verify that the Łojasiewicz exponent of a pair of series is equal to the maximum of the relative exponents. For the exponent at infinity we have only

$$l_\infty(H) \geq \min\{l_\infty(H, X), l_\infty(H, Y)\}. \tag{5}$$

An easy proof of the equality, for the local case, does not transfer to the considered situation if $l_\infty(H) < 0$. However, one can verify that the equality in (5) holds if $l_\infty(H) \geq 0$ or if $H$ is a nondegenerate pair (Corollary 7.4). An example that the inequality can be sharp is $H = (1 + X^4 - Y^2, X^2 - Y)$. By Theorem 5.4 we have $l_\infty(H) = -1$ and by Theorems 5.5 and 5.7 $l_\infty(H, X) = -2$ and $l_\infty(H, Y) = -1$, respectively.

A more convenient version of the relative exponent can be obtained by restriction to a subset of $\mathbf{C}^2$. Let $A \subset \mathbf{C}^2$ be an arbitrary subset. We define $l_\infty(H, A)$ to be the upper bound of the set of all real $\lambda$ such that the inequality (1) holds for $z \in A$ and for sufficiently large $|z|$ ($c > 0$). Obviously, $l_\infty(H, \emptyset) = +\infty$ and $l_\infty(H, \mathbf{C}^2) = l_\infty(H)$. One can easily verify that

$$l_\infty(H, A \cup B) = \min\{l_\infty(H, A), l_\infty(H, B)\}. \tag{6}$$

It is also convenient to consider an exponent related to a variable and to a subset, simultaneously. For $A \subset \mathbf{C}^2$ we define $l_\infty(H, A, X)$ as the upper bound of the set of all real $\lambda$ such that equality (4) holds for $(x, y) \in A$ and for $|x|$ sufficiently large ($c > 0$). Analogously, we define $l_\infty(H, A, Y)$.

A starting point to calculating $l_\infty(H)$ is the following

LEMMA 4.1. *Let fix arbitrary constants $c_1 \geq 1$, $c_2 \geq 1$ and consider the subsets $A = \{|y| \leq c_1|x|\}$ and $B = \{|x| \leq c_2|y|\}$ of $\mathbf{C}^2$. Then*

$$l_\infty(H) = \min\{l_\infty(H, A, X), l_\infty(H, B, Y)\}.$$

Proof. It is easy to verify that $l_\infty(H, A) = l_\infty(H, A, X)$ and $l_\infty(H, B) = l_\infty(H, B, Y)$. Then we use (6) to the equality $A \cup B = \mathbf{C}^2$.

## 5 Laurent-Puiseux series

A convenient tool in a computation of the Łojasiewicz exponent at infinity (in two dimensions) is the classical technique of the Laurent-Puiseux series. These series are a counterpart of the Newton-Puiseux series which are very useful in the local case. In this section we show how Laurent-Puiseux series can be used for calculating the exponents, defined in the previous section.

Denote by $\mathbf{C}((\frac{1}{X}))$ a field of formal Laurent series of variable $X$ with the terms of bounded from above degrees. A *degree* of a series is the maximal power in the expansion, or $-\infty$ for zero series. We say that $p(X) = p_0 X^{\deg p} + \ldots \in \mathbf{C}((\frac{1}{X}))$ ($p_0 \neq 0$) is convergent if it is convergent in a neighbourhood of infinity in $\mathbf{C}$. For every such series we have the inequality

$$c|x|^{\deg p} \leq |p(x)| \leq c'|x|^{\deg p} \qquad (7)$$

for sufficiently large $|x|$, where $0 < c < |p_0| < c'$. Moreover, the differences $c' - |p_0|$ and $|p_0| - c$ can be arbitrary small. For a pair of series $p = (p_1, p_2)$ we define $\deg p = \max\{\deg p_1, \deg p_2\}$. The inequality of the type (7) is also valid for pairs.

The field of formal Laurent-Puiseux series (with quotient powers) is $\mathbf{C}((\frac{1}{X}))^* = \bigcup_{k \geq 1} \mathbf{C}((X^{-\frac{1}{k}}))$. For any non-zero Laurent-Puiseux series $a(X) = a_0 X^{\theta_0} + a_1 X^{\theta_1} + \ldots$ ($\theta_0 > \theta_1 > \ldots$) we put $\deg a(X) = \theta_0$ and $a^+(X) = a_0 X^{\theta_0}$ ($a_0$ will be called *leading coefficient*, and $a^+(X)$ — *leading term*). Moreover, $\deg 0 = -\infty$ and $0^+ = 0$. By the definition, for any $a(X) \in \mathbf{C}((\frac{1}{X}))^*$ there exists a positive integer $d$ such that $a(T^d) \in \mathbf{C}((\frac{1}{T}))$. We say that $a(X)$ is convergent if the corresponding $a(T^d)$ is convergent.

Let $h \in \mathbf{C}[X, Y]$. We say that $a(X) \in \mathbf{C}((\frac{1}{X}))^*$ is a solution of the equation $h(X, Y) = 0$ (with respect to $Y$) if $h(X, a(X)) = 0$ in $\mathbf{C}((\frac{1}{X}))^*$. The minimal positive integer $m$ such that $\frac{\partial^m h}{\partial Y^m}(X, a(X)) \neq 0$ in $\mathbf{C}((\frac{1}{X}))^*$ is called the multiplicity of $a(X)$ and will be denoted by $k(a)$. The polynomial $h$ can be treated as a polynomial of variable $Y$ with coefficients in $\mathbf{C}[X]$. If $h$ is non-zero, then the equation $h(X, Y) = 0$ can be written in the form

$$w_0(X)Y^q + w_1(X)Y^{q-1} + \ldots + w_q(X) = 0, \quad w_i \in \mathbf{C}[X], \ w_0 \neq 0,$$

where $q = \deg_Y f$. In particular, the coefficients $w_i$ are the elements of the field $\mathbf{C}((\frac{1}{X}))^*$. This field is algebraically closed. It comes simply from the algebraic closure of the analogously defined field $\mathbf{C}((X))^*$ ([W], chapter IV, paragraph 3). Hence

$$h(X, Y) = w_0(X) \prod_{a \in \mathcal{H}} (Y - a(X))^{k(a)} , \tag{8}$$

where $\mathcal{H}$ is the set of all the solutions of the equation (we use convention $\prod_\emptyset = 1$). It can be derived, by passing to the local case and by using for example Artin's theorem [A], that all the solutions are convergent. Analogously, we can solve equation $h(X, Y) = 0$ with respect to the variable $X$ in the field $\mathbf{C}((\frac{1}{Y}))^*$.

Further, using factorizations of the form (8), we describe the relative exponents from Lemma 4.1 by using Laurent-Puiseux series (Theorem 5.2 and 5.3). First, we prove a simple

PROPOSITION 5.1. Let $p(T), q_1(T), \ldots, q_k(T)$ be convergent series from $\mathbf{C}((\frac{1}{T}))$ such that $\deg p \leq d$ and $\deg q_j > d$ ($j = 1, \ldots, k$) for a fixed positive integer $d$. Then for any $c > 0$ there exists a neighbourhood $U$ of infinity in $\mathbf{C}$ such that for every $t \in U$ and $|y| \leq c|t|^d$

$$\prod_{i=1}^{k} |y - q_i(t)| \geq \frac{1}{2^k} \prod_{i=1}^{k} |p(t) - q_i(t)| . \tag{9}$$

Proof. From the property of the degree, for fixed $q_i$ there exist $r_i > 0$ such that $|q_i(t)| - c|t|^d \geq \frac{1}{2}|q_i(t) - p(t)|$ for $|t| > r_i$. If additionally $|y| \leq c|t|^d$, then $|y - q_i(t)| \geq |q_i(t)| - |y| \geq |q_i(t)| - c|t|^d \geq \frac{1}{2}|q_i(t) - p(t)|$. Multiplying last inequality for $i = 1, \ldots, k$, we obtain (9) with $U = \{t \in \mathbf{C} : |t| > \max r_i\}$.

Now we are in a position to prove

THEOREM 5.2. Let $H = (f, g)$ be a pair of non-zero polynomials. Denote by $\mathcal{F}_1$ and $\mathcal{G}_1$, respectively, the sets of all the solutions of the equations $f(X, Y) = 0$ and $g(X, Y) = 0$ in $\mathbf{C}((\frac{1}{X}))^*$ with the degrees not greater then 1. Then, there exist a constant $c_1 \geq 1$ such that

$$l_\infty(H, \{|y| \leq c_1|x|\}, X) = \\ \min \{\deg H(X, 0), \inf_{a \in \mathcal{F}_1} \deg g(X, a(X)), \inf_{b \in \mathcal{G}_1} \deg f(X, b(X))\} .$$

Proof. In the proof we use an idea from the Płoski's lemma of the norm of a polynomial mapping ([Pł2], Lemma 3.1).

Denote by $\mathcal{F}$ and $\mathcal{G}$ the sets of all the solutions of the equations $f(X, Y) = 0$ and $g(X, Y) = 0$ in $\mathbf{C}((\frac{1}{X}))^*$, respectively. Let $d$ be a positive integer such that all the considered solution have integer powers after substitution $X = T^d$ (we put d=1 if $\mathcal{F} \cup \mathcal{G} = \emptyset$). Let $c_1 \geq 1$ be a number greater than all the modules of the leading coefficients of all the solutions of the degree 1 from $\mathcal{F} \cup \mathcal{G}$. Put $A = \{|y| \leq c_1|x|\}$. Let $p \in \mathbf{C}((\frac{1}{T}))$ be a convergent series such that $\deg p \leq d$. If $\deg p = d$ then we assume that the module of the leading coefficient of $p$ is less then $c_1$. We claim that

$$l_\infty(H, A, X) \leq \frac{1}{d} \deg H(T^d, p(T)) . \qquad (10)$$

From the property of the degree follows that $(t^d, p(t)) \in A$ for $|t|$ sufficiently large. If $H(T^d, p(T)) = 0$ then (10) holds (both sides are equal $-\infty$). Assume that $H(T^d, p(T)) \neq 0$. Taking an arbitrary $\lambda$ as in the definition of $l_\infty(H, A, X)$, we have $|H(x, y)| \geq c|x|^\lambda$ for $(x, y) \in A$ and $|x|$ sufficiently large ($c > 0$). Then, for sufficiently large $|t|$

$$c'|t|^{\deg H(T^d, p(T))} \geq |H(t^d, p(t))| \geq c|t^d|^\lambda .$$

This means that $\deg H(T^d, p(T)) \geq \lambda d$, what gives (10), considering an arbitrary choice of $\lambda$. Now, the inequality ($\leq$) in the theorem can be obtained by taking as $p(T)$ the zero series or any solution from $\mathcal{F}_1 \cup \mathcal{G}_1$ after substitution $X = T^d$.

Before we prove the opposite inequality, let us make the following observation. Let $\{p(T)\}$ be a finite family, where $p(T) \in \mathbf{C}((\frac{1}{T}))$ is a convergent non-zero series, or $p(T) \in \mathbf{C}((\frac{1}{T})) \times \mathbf{C}((\frac{1}{T}))$ is a non-zero pair of convergent series. Notice, that there exist a neighbourhood $V$ of infinity in $\mathbf{C}$ and a constant $c > 0$ such that for every $t \in V$ and for every $p$ from this family

$$|p(t)| \geq c|t|^{\deg p} .$$

Now, consider the family which is composed of the pair $H(T^d, 0)$ and all the series $f(T^d, b(T))$ for $b \in \mathcal{G}_1$, and $g(T^d, a(T))$ for $a \in \mathcal{F}_1$. If some element of the family is a zero-element, then the both sides of the equality in the statement of the theorem are equal to $-\infty$, and the theorem is true. Therefore, we can assume that all the elements of the family are non-zero. Let $V$ be a neighbourhood of infinity in $\mathbf{C}$ and let $c > 0$ be a constant constructed for this family as above. Let $u_0(X)$ be the coefficient at $Y^{\deg_Y f}$ if we treat $f$ as a polynomial of variable $Y$ with coefficients in $\mathbf{C}[X]$. By (8), we can write $f = f_1 f_2$, where $f_1$ is the product of $(Y - a(X))^{k(a)}$ for $a \in \mathcal{F}_1$ ($\deg a \leq 1$), and $f_2$ is the product of such factors for $a \in \mathcal{F} \setminus \mathcal{F}_1$ ($\deg a > 1$) and $u_0(X)$. We consider an analogous factorization for $g = g_1 g_2$, where where $g_1$ is the product of $(Y - b(X))^{k(a)}$ for $b \in \mathcal{G}_1$ ($\deg b \leq 1$), and $g_2$ is the product of such factors for $b \in \mathcal{G} \setminus \mathcal{G}_1$ ($\deg b > 1$) and $v_0 \in \mathbf{C}[X]$ analogously defined for $g$. Applying Proposition 5.1 we can find a neighbourhood $U$ of infinity in $\mathbf{C}$ such that for $t \in U$ and $|y| \leq c_1|t|^d$

$$|f_2(t^d, y)| \geq 2^{-\deg_Y f_2} |f_2(t^d, b(t^d))| \qquad (11)$$

for every $b \in \mathcal{G}_1 \cup \{0\}$ and

$$|g_2(t^d, y)| \geq 2^{-\deg_Y g_2}|g_2(t^d, a(t^d))| \tag{12}$$

for every $a \in \mathcal{F}_1 \cup \{0\}$. Now, take an arbitrary pair $(t, y)$ such that $t \in U \cap V$ and $|y| \leq c_1|t|^d$. Denote by $m$ the minimum on the right side of the equality in the statement of the theorem. If $\mathcal{F}_1 = \mathcal{G}_1 = \emptyset$, then $f = f_2$, $g = g_2$ and $\deg_Y f = \deg_Y f_2$, $\deg_Y g = \deg_Y g_2$. Applying (11) with $b = 0$ and (12) with $a = 0$ we have

$$\begin{aligned} |H(t^d, y)| &= \max\{|f(t^d, y)|, |g(t^d, y)|\} \geq \\ &\geq 2^{-\max\{\deg_Y f, \deg_Y g\}} \max\{|f(t^d, 0)|, |g(t^d, 0)\} \geq \\ &\geq 2^{-\deg_Y H} c|t^d|^{\deg H(X,0)} \geq 2^{-\deg_Y H} c|t^d|^m \ . \end{aligned}$$

If $\mathcal{F}_1 = \emptyset$ and $\mathcal{G}_1 \neq \emptyset$, then $f = f_2$ and $\deg_Y f = \deg_Y f_2$. Let $b \in \mathcal{G}_1$. Using (11), we have

$$\begin{aligned} |H(t^d, y)| &\geq |f(t^d, y)| \geq 2^{-\deg_Y f}|f(t^d, b(t^d))| \geq \\ &\geq 2^{-\deg_Y f} c|t|^{\deg f(T^d, b(T^d))} \geq 2^{-\deg_Y H} c|t^d|^m \ . \end{aligned}$$

The case $\mathcal{F}_1 \neq \emptyset$ and $\mathcal{G}_1 = \emptyset$ is analogous. If both the sets $\mathcal{F}_1$ and $\mathcal{G}_1$ are non-empty, then we consider the finite subsets of complex numbers $\mathcal{F}_1^t = \{a(t^d) : a \in \mathcal{F}_1\}$ and $\mathcal{G}_1^t = \{b(t^d) : b \in \mathcal{G}_1\}$, for earlier fixed $t$. Suppose that $\text{dist}(\mathcal{F}_1^t, y) \geq \text{dist}(\mathcal{G}_1^t, y)$. Let $\tilde{b}(t^d) \in \mathcal{G}_1^t$ be a point realizing the distance $\text{dist}(y, \mathcal{G}_1^t)$ for $\tilde{b} \in \mathcal{G}_1$. For every $a \in \mathcal{F}_1$ we have $|y - a(t^d)| \geq \text{dist}(y, \mathcal{F}_1^t) \geq \text{dist}(y, \mathcal{G}_1^t) = |y - \tilde{b}(t^d)|$, hence $|y - a(t^d)| \geq \frac{1}{2}|y - a(t^d)| + \frac{1}{2}|y - \tilde{b}(t^d)| \geq \frac{1}{2}|\tilde{b}(t^d) - a(t^d)|$, and therefore

$$|f_1(t^d, y)| \geq 2^{-\deg_Y f_1}|f_1(t^d, \tilde{b}(t^d))| \ ,$$

what in addition to (11) gives $|f(t^d, y)| \geq 2^{-\deg_Y f}|f(t^d, \tilde{b}(t^d))|$. Hence

$$|H(t^d, y)| \geq |f(t^d, y)| \geq 2^{-\deg_Y f} c|t|^{\deg f(T^d, \tilde{b}(T^d))} \geq 2^{-\deg H} c|t^d|^m \ .$$

The case $\text{dist}(y, \mathcal{F}_1^t) \leq \text{dist}(y, \mathcal{G}_1^t)$ can be verified analogously. This ends the proof of the theorem.

By symmetry we obtain

THEOREM 5.3. *Let $H = (f, g)$ be a pair of non-zero polynomials. Denote by $\mathcal{F}_1$ and $\mathcal{G}_1$, respectively, the sets of all the solutions of the equations $f(X, Y) = 0$ and $g(X, Y) = 0$ in $\mathbf{C}((\frac{1}{Y}))^*$ with the degrees not greater then 1. Then, there exist a constant $c_2 \geq 1$ such that*

$$l_\infty(H, \{|x| \leq c_2|y|\}, Y) = \\ \min\{\deg H(0, Y), \inf_{a \in \mathcal{F}_1} \deg g(a(Y), Y), \inf_{b \in \mathcal{G}_1} \deg f(b(Y), Y)\} \ .$$

Joining Theorems 5.2 and 5.3 with Lemma 4.1 we obtain a formula for calculating $l_\infty(H)$ in an arbitrary situation

THEOREM 5.4. *Let $H = (f, g)$ be a pair of non-zero polynomials. Denote by $\mathcal{F}'_1$ and $\mathcal{G}'_1$, respectively, the sets of all the solutions of the equations $f(X, Y) = 0$ and $g(X, Y) = 0$ in $\mathbf{C}((\frac{1}{X}))^*$ with the degrees not greater than 1, and by $\mathcal{F}''_1$ and $\mathcal{G}''_1$ the sets of all the solutions of the same equations in $\mathbf{C}((\frac{1}{Y}))^*$ with the degrees also not greater than 1. Then the exponent $l_\infty(H)$ is equal to the minimum of the following six quantities*

$$\deg H(X, 0), \inf_{a \in \mathcal{F}'_1} \deg g(X, a(X)), \inf_{b \in \mathcal{G}'_1} \deg f(X, b(X)),$$

$$\deg H(0, Y), \inf_{a \in \mathcal{F}''_1} \deg g(a(Y), Y), \inf_{b \in \mathcal{G}''_1} \deg f(b(Y), Y).$$

This type result was obtained by Chądzyński & Krasiński in [ChK1].

The following theorem unables us to calculate $l_\infty(H, X)$.

THEOREM 5.5. *For a pair $H = (f, g)$ of non-zero polynomials*

$$l_\infty(H, X) = \min \left\{ \deg H(X, 0), \inf_{a \in \mathcal{F}} \deg g(X, a(X)), \inf_{b \in \mathcal{G}} \deg f(X, b(X)) \right\},$$

*where $\mathcal{F}$ and $\mathcal{G}$ are the sets of all the solutions of the equations $f(X, Y) = 0$ and $g(X, Y) = 0$ in $\mathbf{C}((\frac{1}{X}))^*$, respectively.*

Proof. Let $\mathcal{F}, \mathcal{G}$ be the sets of all the solutions of the equations $f(X, Y) = 0$, $g(X, Y) = 0$ in $\mathbf{C}((\frac{1}{X}))^*$, respectively. Let $d$ be a positive integer such that all the sollutions have integer powers after substitution $X = T^d$ ($d = 1$ in $\mathcal{F} \cup \mathcal{G} = \emptyset$). The inequality "$\leq$", in the theorem, can be obtained similarily as in the proof of Theorem 5.2, by using the estimate $l_\infty(H, X) \leq \frac{1}{d} \deg H(T^d, p(T))$ for $p(T) \in \mathbf{C}((\frac{1}{T}))$. The opposite inequality can be obtained by using the idea which comes from Płoski [Pł2]. We use so called "lemma of the norm of a polynomial mapping".

LEMMA 5.6 ([Pł2] – Lemma 3.1, [L] – Lemma 4.7).
*If $w = (u, v)$ is a pair of non-zero polynomials $u(Y), v(Y) \in \mathbf{C}[Y]$, then for every $y \in \mathbf{C}$*

$$\max \{|u(y)|, |v(y)|\} \geq 2^{-n} \min \left\{ |w(0)|, \inf_{\eta \in u^{-1}(0)} |v(\eta)|, \inf_{\eta \in v^{-1}(0)} |u(\eta)| \right\},$$

*where $n = \max \{\deg u, \deg v\}$.*

Now, cosider the family of series $H(T^d, 0)$, $f(T^d, b(T^d))$ for $b \in \mathcal{G}$, and $g(T^d, a(T^d))$ for $a \in \mathcal{F}$. As in the proof of Theorem 5.2 we can choose $c > 0$ and a neighbourhood $V$ of infinity in $\mathbf{C}^2$ such that $|p(t)| \geq c|t|^{\deg p}$ for every $t \in V$ and for every member $p$ of the family. Now, fix $t \in V$ and $y \in \mathbf{C}$. Let $m$ be the minimum from the statement of the theorem. Using factorizations for $f$ and $g$, of the form (8), and Lemma 5.6, we obtain

$$\max\{|f(t^d, y)|, |g(t^d, y)|\} \geq$$
$$\geq 2^{-\deg_Y H} \min\left\{|H(t^d, 0)|, \inf_{a \in \mathcal{F}} |g(t^d, a(t^d))|, \inf_{b \in \mathcal{G}} |f(t^d, b(t^d))|\right\}$$
$$\geq c \cdot 2^{-\deg_Y H} |t^d|^m .$$

This ends the proof of the theorem.

By symmetry we obtain

THEOREM 5.7. *For a pair $H = (f, g)$ of non-zero polynomials*

$$l_\infty(H, Y) = \min\left\{\deg H(0, Y), \inf_{a \in \mathcal{F}} \deg g(a(Y), Y), \inf_{b \in \mathcal{G}} \deg f(b(Y), Y)\right\},$$

*where $\mathcal{F}$ and $\mathcal{G}$ are the sets of all the solutions of the equations $f(X, Y) = 0$ and $g(X, Y) = 0$ in $\mathbf{C}((\frac{1}{Y}))^*$, respectively.*

# 6 Relations with Newton diagram

In this section we describe connections between the degrees of the Laurent-Puiseux solutions of the equation $h(X, Y) = 0$ and the Newton diagram $\Delta_h$ for a non-zero polynomial $h$. The results are classical. The local case is well described (see, for example [W], [BK], [Pł2]). We focus our attention on solving the equation with respect to $Y$ in $\mathbf{C}((\frac{1}{X}))^*$. The facts concerning solutions with respect to $X$ in $\mathbf{C}((\frac{1}{Y}))^*$ can be obtained by symmetry. We start from the following

PROPOSITION 6.1. *Let $h \in \mathbf{C}[X, Y]$ be non-zero polynomial and let $S \in \mathcal{N}_h^{(\mathbf{r})}$. Then:*

(a) *There exist a factorization of $\mathrm{in}(h, S)$ in $\mathbf{C}((\frac{1}{X}))^*[Y]$ of the form*

$$\mathrm{in}(h, S) = \epsilon X^\zeta Y^\vartheta \prod_{i=1}^{|S_2|} \left(Y - a_i X^{\frac{|S_1|}{|S_2|} \sigma(S)}\right) ,$$

*where $\zeta = \sigma(S) \min\{\alpha\sigma(S) : (\alpha, \beta) \in S\}$, $\vartheta = \min\{\beta : (\alpha, \beta) \in S\}$ and $\epsilon, a_1, \ldots, a_{|S_2|}$ are non-zero complex numbers.*

(b) *If $\theta$ is a rational number and $c$ is a non-zero complex number such that $\mathrm{in}(h, S)(X, cX^\theta) = 0$, then $\theta = \frac{|S_1|}{|S_2|} \sigma(S)$ and $c$ is one of the numbers $a_1, \ldots, a_{|S_2|}$ described in (a).*

Proof. Notice that (b) follows immediately from (a). To prove (a) consider a factorization of $\mathrm{in}(h, S)$ of the form

$$\mathrm{in}(h, S) = X^\zeta Y^\vartheta \prod_{i=1}^{d} (u_i X^{\xi \sigma(S)} + v_i Y^\eta) ,$$

where $d = \mathrm{GCD}\{|S_1|, |S_2|\}$, $\xi = |S_1|/d$, $\eta = |S_2|/d$ and $u_i$, $v_i$ are non-zero complex numbers. Consider a factor of the form $uX^{\xi\sigma(S)} + vY^\eta$ and a new variable $X' = X^{\frac{\xi}{\eta}\sigma(S)} = X^{\frac{|S_1|}{|S_2|}\sigma(S)}$. Then $uX^{\xi\sigma(S)} + vY^\eta = uX'^\eta + vY^\eta = v(Y - \epsilon_1 X') \ldots (Y - \epsilon_\eta X')$, where $\epsilon_1, \ldots, \epsilon_\eta$ are the complex roots of $-\frac{u}{v}$ of the degree $\eta$. This ends the proof.

Now, fix an arbitrary rational $\theta$ and consider a linear form $(\alpha, \beta) \mapsto \alpha + \theta\beta$. Since $\operatorname{supp} h$ is nonempty the subset of $\operatorname{supp} h$, where the form attains its maximum, is also non-empty. Let $m = \max\{\alpha + \theta\beta : (\alpha, \beta) \in \operatorname{supp} h\}$. Consider a polynomial

$$\sum_{\alpha+\theta\beta=m} h_{\alpha\beta} X^\alpha Y^\beta . \tag{13}$$

The equation $\alpha + \theta\beta = m$ describes the line supporting $\Delta_h$. The polynomial (13) is a monomial, if the line meets the diagram at one point, or a quasi homogeneous form $\operatorname{in}(h, S)$, if the line meets the diagram along a segment $S \in \mathcal{N}_h^{(\mathbf{r})}$. The second possibility occurs exactly when the maximum of the considered linear form is attained at more than one point of $\operatorname{supp} h$. One can verify that it happens if and only if $\theta$ is one of the numbers $\frac{|S_1|}{|S_2|}\sigma(S)$, $S \in \mathcal{N}_h^{(\mathbf{r})}$. One can treat the maximum $m$ as the maximal possible degree of the substitution $h(X, a(X))$, where $a \in \mathbf{C}((\frac{1}{X}))^*$ and $\deg a = \theta$. If we put $H_{\alpha\beta} = h_{\alpha\beta} X^\alpha a(X)^\beta$ for any $h_{\alpha\beta} \neq 0$, then $h(X, a(X)) = \sum H_{\alpha\beta}(X)$, $H_{\alpha\beta}^+ = h_{\alpha\beta} X^\alpha a^+(X)^\beta$ and $\deg H_{\alpha\beta} = \alpha + \theta\beta$. Using standard properties of the degree we obtain

$$\deg h(X, a(X)) \leq \max\{\deg H_{\alpha\beta}\} = m . \tag{14}$$

This observation gives a motivation for a definition. We say that a substitution $h(X, a(X))$ is *generic* if the equality at (14) holds, and conversely, *non-generic* if the strong inequality holds.

PROPOSITION 6.2. *If $h \in \mathbf{C}[X, Y]$ is a non-zero polynomial and $a(X)$ is a non-zero Laurent-Puiseux series, then the substitution $h(X, a(X))$ is non-generic if and only if there exist a segment $S \in \mathcal{N}_h^{(\mathbf{r})}$ such that $\operatorname{in}(h, S)(X, a^+(X)) = 0$.*

Proof. We use notation as earlier. Notice that

$$h(X, a(X)) = \sum_{\alpha+\theta\beta=m} H_{\alpha\beta}^+ + \{\text{terms of lower degree}\} . \tag{15}$$

If the substitution is non-generic, then the first component of the sum (15) vanish. It means that the number of terms with the maximal degree $m$ is greater than one (to obtain a reduction). In such a case the sum (13) is equal to $\operatorname{in}(h, S)$ for a segment $S \in \mathcal{N}_h^{(\mathbf{r})}$, and then

$$0 = \sum_{\alpha+\theta\beta=m} H_{\alpha\beta}^+ = \sum_{\alpha+\theta\beta=m} h_{\alpha\beta} X^\alpha a^+(X)^\beta = \operatorname{in}(h, S)(X, a^+(X)) ,$$

what ends the proof of ($\Rightarrow$). To prove the opposite implication assume that there exists a segment $S \in \mathcal{N}_h^{(\mathbf{r})}$ such that $\mathrm{in}(h,S)(X, a^+(X)) = 0$. By Proposition 4.8(b) $\deg a(X) = \frac{|S_1|}{|S_2|}\sigma(S)$. Let $\theta = \frac{|S_1|}{|S_2|}\sigma(S)$, and let $m$ be the maximum of $\alpha + \theta\beta$ over $(\alpha, \beta) \in \mathrm{supp}\, h$, as earlier. The sum (13) is equal to $\mathrm{in}(h,S)$, in this case. The equality $\mathrm{in}(X, a^+(X)) = 0$ means that the first component in the right side of (15) vanish, therefore the substitution $h(X, a(X))$ is non-generic.

The following theorem is a counterpart of the classical local result (for example see [Pł2], [BK]).

THEOREM 6.3. *Let $h \in \mathbf{C}[X,Y]$ be a non-zero polynomial. Then:*

(a) *If $a(X)$ is a non-zero solution of the equation $h(X,Y) = 0$ in $\mathbf{C}((\frac{1}{X}))^*$ then there exists a segment $S \in \mathcal{N}_h^{(\mathbf{r})}$ such that $\deg a = \frac{|S_1|}{|S_2|}\sigma(S)$ and $\mathrm{in}(h,S)(X, a^+(X)) = 0$.*

(b) *For every segment $S \in \mathcal{N}_h^{(\mathbf{r})}$ there exist exactly $|S_2|$ non-zero solutions in $\mathbf{C}((\frac{1}{X}))^*$ of the equation $h(X,Y) = 0$ of the degree $\theta = \frac{|S_1|}{|S_2|}\sigma(S)$, counting with multiplicities. Moreover, if $(Y - cX^\theta)^k$ is a factor of $\mathrm{in}(h,S)$ in $\mathbf{C}((\frac{1}{X}))^*[Y]$, then there exist exactly $k$ solutions, counting with multiplicities, with leading term $cX^\theta$.*

Proof. The part (a) of the theorem follows immediately from Proposition 6.2. To prove (b) fix a segment $S \in \mathcal{N}_h^{(\mathbf{r})}$. For any non-zero $f \in \mathbf{C}((\frac{1}{X}))^*[Y]$ we can define weighted degree

$$\deg{}^*f = \max\left\{\alpha + \beta\frac{|S_1|}{|S_2|}\sigma(S) : (\alpha, \beta) \in \mathrm{supp}\, f\right\}.$$

This number is well defined, because it is the maximum of the bounded from above set of rationals with the bounded denominators. Define $f^* = \{$*sum of monomials with the maximal weighted degree*$\}$. We put $deg^*0 = -\infty$ and $0^* = 0$. The above-defined weighted degree has standard properties: $\deg{}^*(fg) = \deg{}^*f + \deg{}^*g$ and $(fg)^* = f^*g^*$. Notice that $h^* = \mathrm{in}(h,S)$. Now, denote by $\mathcal{H}$ the set of all the solutions of $h(X,Y) = 0$ in $\mathbf{C}((\frac{1}{X}))^*$ and consider a factorization (8) for $h$. By the property of the weighted degree we have

$$h^* = \mathrm{in}(h,S) = w_0(X)^* \prod_{a \in \mathcal{H}} [(Y - a(X))^*]^{k(a)},$$

where $w_0(X)^* = w_0(X)^+$ and

$$(Y - a(X))^* = \begin{cases} Y & \text{if } \deg a(X) < \frac{|S_1|}{|S_2|}\sigma(S), \\ Y - a^+(X) & \text{if } \deg a(X) = \frac{|S_1|}{|S_2|}\sigma(S), \\ -a^+(X) & \text{if } \deg a(X) > \frac{|S_1|}{|S_2|}\sigma(S). \end{cases}$$

Now the thesis (b) follows immediately from the unique factorization of $\mathbf{C}((\frac{1}{X}))^*[Y]$ and from Proposition 6.1(a).

By symmetry we can obtain the results concerning solutions of the equation $h(X,Y) = 0$ in $\mathbf{C}((\frac{1}{Y}))^*$. From Proposition 6.2 we have

PROPOSITION 6.4. *If $h \in \mathbf{C}[X,Y]$ is a non-zero polynomial and $a(Y)$ is a non-zero Laurent-Puiseux series, then the substitution $h(a(Y), Y)$ is non-generic if and only if there exist a segment $S \in \mathcal{N}_h^{(\mathbf{t})}$ such that $\mathrm{in}(h, S)(a^+(Y), Y) = 0$.*

From Theorem 6.3 we obtain

THEOREM 6.5. *Let $h \in \mathbf{C}[X,Y]$ be a non-zero polynomial. Then:*

(a) *If $a(Y)$ is a non-zero solution of the equation $h(X,Y) = 0$ in $\mathbf{C}((\frac{1}{Y}))^*$ then there exists a segment $S \in \mathcal{N}_h^{(\mathbf{t})}$ such that $\deg a = \frac{|S_2|}{|S_1|}\sigma(S)$ and $\mathrm{in}(h, S)(a^+(Y), Y) = 0$.*

(b) *For every segment $S \in \mathcal{N}_h^{(\mathbf{t})}$ there exist exactly $|S_1|$ non-zero solutions in $\mathbf{C}((\frac{1}{Y}))^*$ of the equation $h(X,Y) = 0$ of the degree $\theta = \frac{|S_2|}{|S_1|}\sigma(S)$, counting with the multiplicities. Moreover, if $(X - cY^\theta)^k$ is a factor of $\mathrm{in}(h, S)$ in $\mathbf{C}((\frac{1}{Y}))^*[Y]$, then there exist exactly $k$ solutions, counting with multiplicities, with leading term $cY^\theta$.*

# 7 Proofs of auxiliary results

In this section we prove Theorems 3.1 and 3.2 by using Theorem 5.4 and the facts from the previous section. Let $h \in \mathbf{C}[X,Y]$ be a non-zero polynomial. Notice, that Theorem 6.3 determines a correspondence between non-zero solutions of $h(X,Y) = 0$ in $\mathbf{C}((\frac{1}{X}))^*$ and the segments of $\mathcal{N}_h^{(\mathbf{r})}$. To each solution $a(X)$ we can assign a segment $S \in \mathcal{N}_h^{(\mathbf{r})}$ such that $\deg a(X) = \frac{|S_1|}{|S_2|}\sigma(S)$. Analogously, Theorem 6.5 determines similar correspondence between non-zero solutions of the equation in $\mathbf{C}((\frac{1}{Y}))^*$ and the segments of $\mathcal{N}_h^{(\mathbf{t})}$.

To prove Theorem 3.1 consider a pair $H = (f, g)$ of non-zero polynomials. Let $a(X) \in \mathbf{C}((\frac{1}{X}))^*$ be a non-zero solution of the equation $f(X,Y) = 0$ and let $S \in \mathcal{N}_f^{(\mathbf{r})}$ corresponds to this solution ($\deg a(X) = \frac{|S_1|}{|S_2|}\sigma(S)$). Since the number $\alpha(S, \Delta_g)$, defined by (2), is exactly equal to the maximal possible degree of the substitution $g(X, a(X))$, we have the inequality

$$\deg g(X, a(X)) \leq \alpha(S, \Delta_g) . \qquad (16)$$

Now, we show that the equality holds if the pair is nondegenerate. To obtain a contradiction suppose that the inequality is strong. Then the substitution is non-generic and by Proposition 6.2 there exists a segment $T \in \mathcal{N}_g^{(\mathbf{r})}$ such that $\text{in}(g, T)(X, a^+(X)) = 0$. Moreover, by Proposition 6.1(b) $\frac{|T_1|}{|T_2|}\sigma(T) = \deg a(X)$, so $T$ and $S$ are parallel. By Theorem 6.5(a) we have also $\text{in}(f, S)(X, a^+(X)) = 0$. Hence, the system of equations $\text{in}(f, S) = \text{in}(g, T) = 0$ has a solution in $(\mathbf{C} \setminus \{0\}) \times (\mathbf{C} \setminus \{0\})$, what means degeneracy of the pair.

Now, consider the set $\mathcal{F}_1'$ of all the solutions of $f(X, Y) = 0$ in $\mathbf{C}((\frac{1}{X}))^*$ of the degrees less or equal to 1. Let $\mathcal{N}_f'$ be the set of $S \in \mathcal{N}_f^{(\mathbf{r})}$ such that $\frac{|S_1|}{|S_2|}\sigma(S) \leq 1$. For every non-zero solution in $\mathcal{F}_1'$ there exists a segment of $\mathcal{N}_f'$ such that the declivity of the segment is equal to the degree of the solution. From (16) we obtain

$$\inf_{a \in \mathcal{F}_1'} \deg g(X, a(X)) \leq \inf_{a \in \mathcal{F}_1' \setminus \{0\}} \deg g(X, a(X)) \leq \inf_{S \in \mathcal{N}_f'} \alpha(S, \Delta_g) \, .$$

Let $\mathcal{F}_1''$ be the set of all the solutions of $f(X, Y) = 0$ in $\mathbf{C}((\frac{1}{Y}))^*$ of the degrees less or equal to 1. There is a correspondence between non-zero solutions of $\mathcal{F}_1''$ and a polygon $\mathcal{N}_f''$ which consists of $S \in \mathcal{N}_f^{(\mathbf{t})}$, $\frac{|S_2|}{|S_1|}\sigma(S) \leq 1$. We define, analogously, the sets of solutions $\mathcal{G}_1'$ and $\mathcal{G}_1''$ for the equation $g(X, Y) = 0$ and the polygons $\mathcal{N}_g'$ and $\mathcal{N}_g''$. By repeating earlier considerations, to the rest infimae from the statement of Theorem 5.4, we obtain $l_\infty(H)$ to be less or equal to the minimum of the six quantities

$$\deg H(X, 0), \quad \inf_{S \in \mathcal{N}_f'} \alpha(S, \Delta_g), \quad \inf_{T \in \mathcal{N}_g'} \alpha(T, \Delta_f), \tag{17}$$

$$\deg H(0, Y), \quad \inf_{S \in \mathcal{N}_f''} \beta(S, \Delta_g), \quad \inf_{T \in \mathcal{N}_g''} \beta(T, \Delta_f). \tag{18}$$

Now, we prove the equality in the case of nondegeneracy. Denote by $m$ the minimum of the above six quantities. It is sufficient to show that each of the six numbers from the statement of Theorem 5.4 is greater or equal to $m$. For $\deg H(X, 0)$ and $\deg H(0, Y)$ it is obvious. Using equality in (16), which is proved for nondegenerate pairs, we obtain

$$\inf_{a \in \mathcal{F}_1'} \deg g(X, a(X)) = \begin{cases} \inf_{S \in \mathcal{N}_f'} \alpha(S, \Delta_g) & \text{if } 0 \notin \mathcal{F}_1' , \\ \min \{ \inf_{S \in \mathcal{N}_f'} \alpha(S, \Delta_g), \deg g(X, 0) \} & \text{if } 0 \in \mathcal{F}_1' . \end{cases}$$

In the above formula we need the fact that for each segment from $\mathcal{N}_f'$ there exists a corresponding solution from $\mathcal{F}_1'$ (Theorem 6.3(b)). If $0 \in \mathcal{F}_1'$, then $\deg g(X, 0) = \deg H(X, 0)$. Hence, in both the cases the left infimum in the formula is greater or equal to $m$. Applying above consideration, to the rest infimae from the statement of Theorem 5.4, we end the proof of the desired equality. In order to obtain a complete proof of Theorems 3.1 and 3.2 it is enough to prove the following

PROPOSITION 7.1. *Let $H = (f,g)$ be a pair of non-zero polynomials. Denote by $m_1$ the minimum of the six quantities from the statement of Theorem 3.1, and by $m_2$ the minimum of the six quantities in (17) and (18). We claim that $m_1 = m_2$.*

Proof. Considering inclusions of the type $\mathcal{N}' \subset \mathcal{N}^{(\mathbf{r})}$ and $\mathcal{N}'' \subset \mathcal{N}^{(\mathbf{t})}$ we have $m_1 \leq m_2$. In order to prove the opposite inequality let us write

$$\inf_{S \in \mathcal{N}_f^{(\mathbf{r})}} \alpha(S, \Delta_g) = \min \left\{ \inf_{S \in \mathcal{N}_f'} \alpha(S, \Delta_g) \,, \; \inf_{S \in \mathcal{N}_f^{(\mathbf{r})} \setminus \mathcal{N}_f'} \alpha(S, \Delta_g) \right\} .$$

Assume that $S \in \mathcal{N}_f^{(\mathbf{r})} \setminus \mathcal{N}_f'$. Then $\frac{|S_1|}{|S_2|} \sigma(S) > 1$, and in particular $\sigma(S) = 1$. Hence

$$\alpha + \beta \frac{|S_1|}{|S_2|} \sigma(S) = \frac{|S_1|}{|S_2|} \sigma(S) \left( \alpha \frac{|S_2|}{|S_1|} \sigma(S) + \beta \right) \geq \alpha \frac{|S_2|}{|S_1|} \sigma(S) + \beta .$$

Considering (2) and (3) we obtain $\alpha(S, \Delta_g) \geq \beta(S, \Delta_g)$. By the obvious inclusion $\mathcal{N}_f^{(\mathbf{r})} \setminus \mathcal{N}_f' \subset \mathcal{N}_f''$ we have

$$\inf_{S \in \mathcal{N}_f^{(\mathbf{r})} \setminus \mathcal{N}_f'} \alpha(S, \Delta_g) \geq \inf_{S \in \mathcal{N}_f^{(\mathbf{r})} \setminus \mathcal{N}_f'} \beta(S, \Delta_g) \geq \inf_{S \in \mathcal{N}_f''} \beta(S, \Delta_g) .$$

Finally

$$\inf_{S \in \mathcal{N}_f^{(\mathbf{r})}} \alpha(S, \Delta_g) \geq \min \left\{ \inf_{S \in \mathcal{N}_f'} \alpha(S, \Delta_g) \,, \; \inf_{S \in \mathcal{N}_f''} \beta(S, \Delta_g) \right\} .$$

Applying analogous considerations to the rest three infimae we obtain desired inequality $m_1 \geq m_2$. This ends the proof of the proposition and the proofs of Theorems 3.1 and 3.2.

From Theorem 5.5 and (16) we obtain

THEOREM 7.2. *For a pair $H = (f,g)$ of non-zero polynomials*

$$l_\infty(H, X) \leq \min \left\{ \deg H(X, 0), \inf_{S \in \mathcal{N}_f^{(\mathbf{r})}} \alpha(S, \Delta_g), \inf_{T \in \mathcal{N}_g^{(\mathbf{r})}} \alpha(T, \Delta_f) \right\} ,$$

*with the equality for the nondegenerate pair.*

By symmetry, we have

THEOREM 7.3. *For a pair $H = (f,g)$ of non-zero polynomials*

$$l_\infty(H, Y) \leq \min \left\{ \deg H(0, Y), \inf_{S \in \mathcal{N}_f^{(\mathbf{t})}} \beta(S, \Delta_g), \inf_{T \in \mathcal{N}_g^{(\mathbf{t})}} \beta(T, \Delta_f) \right\} ,$$

*with the equality for the nondegenerate pair.*

A simple consequence of Theorem 3.2 and the two above theorems is

COROLLARY 7.4. *If a pair $H = (f, g)$ of non-zero polynomials is nondegenerate then*
$$l_\infty(H) = \min\{l_\infty(H, X),\, l_\infty(H, Y)\}\ .$$

# 8 Newton diagrams of derivatives

The theorems presented in the previous section gives us the estimation of $l_\infty(\operatorname{grad} h)$ for a polynomial $h(X, Y)$ based on the information from the Newton diagrams of the gradient's components. Our aim is to give an estimation of $l_\infty(\operatorname{grad} h)$ by using an information from the Newton diagram of $h$, only. We can do this by describing a structure of the right and top Newton polygons of the derivatives $\frac{\partial h}{\partial X}$ and $\frac{\partial h}{\partial Y}$. This is the main goal of this section. We concentrate our attention on describing a structure of the right polygons. The analogous description for the top polygons can be obtained by symmetry, but will not be expressed explicidely.

We start from the description of the right Newton polygon of $\frac{\partial h}{\partial Y}$. It is easy to observe that $\operatorname{supp}\frac{\partial h}{\partial Y}$ is the image of the set $\operatorname{supp} h \setminus \{\beta = 0\}$ in the translation $(\alpha, \beta) \mapsto (\alpha, \beta - 1)$. We say that a segment $T$ of $\mathcal{N}^{(\mathbf{r})}_{\frac{\partial h}{\partial Y}}$ is a *standard* one if it is the image of a segment $S$ of $\mathcal{N}^{(\mathbf{r})}_h$ in this translation. We will use notation $T = S - (0, 1)$. We have

$$\operatorname{in}\left(\frac{\partial h}{\partial Y}, T\right) = \frac{\partial}{\partial Y}\operatorname{in}(h, S)\ , \tag{19}$$

in this case. Obviously, non-standard segments can also exist. Now, we are going to describe the structure of standard and non-standard segments of the considered polygon. We put the right vertices of $\Delta_h$ in order with respect to its ordinates. The first one is the nearest to the horizontal axis.

Let us note the following simple facts. The right Newton polygon of a polynomial is nonempty if and only if its support contains at least two points with the different ordinates. We obtain from this that if the right polygon of $\frac{\partial h}{\partial Y}$ is nonempty than the right polygon of $h$ is also nonempty. So, a necessary condition for $\mathcal{N}^{(\mathbf{r})}_{\frac{\partial h}{\partial Y}}$ to be nonempty is $\mathcal{N}^{(\mathbf{r})}_h \neq \emptyset$. For any right vertex of $\Delta_h$ with a positive ordinate, the image of this vertex in the translation $(\alpha, \beta) \mapsto (\alpha, \beta - 1)$ is a right vertex of $\Delta_{\frac{\partial h}{\partial Y}}$. Similarly, if a segment $S \in \mathcal{N}^{(\mathbf{r})}_h$ does not touch the horizontal axis, then $T = S - (0, 1)$ is a standard segment of $\mathcal{N}^{(\mathbf{r})}_{\frac{\partial h}{\partial Y}}$. In particular, if the first right vertex of $\Delta_h$ does not lay on the horizontal axis, then all the segments of $\mathcal{N}^{(\mathbf{r})}_{\frac{\partial h}{\partial Y}}$ are standard. So, a necessary condition for $\mathcal{N}^{(\mathbf{r})}_{\frac{\partial h}{\partial Y}}$ to has non-standard segments is: the first right vertex of $\Delta_h$ has to lay on the horizontal axis. Consider a right vertex $(\mu, \nu)$ of $\Delta_h$ with the minimal positive ordinate. From the above considerations follows that each segment of $\mathcal{N}^{(\mathbf{r})}_{\frac{\partial h}{\partial Y}}$ over the vertex $(\mu, \nu - 1)$ is standard, and each segment below this wertex is non-standard. So, non-standard segment of $\mathcal{N}^{(\mathbf{r})}_{\frac{\partial h}{\partial Y}}$ can only exist in the strip $\{0 \leq \beta \leq \nu - 1\}$. Now, we give a sufficient condition for $\mathcal{N}^{(\mathbf{r})}_{\frac{\partial h}{\partial Y}}$ to have non-standard segments.

PROPOSITION 8.1. *Suppose that the polygon $\mathcal{N}^{(\mathbf{r})}_h$ is nonempty and the first right vertex of $\Delta_h$ lies on the horizontal axis. Then, this vertex is the lower end of the first segment $F$ of the polygon and the upper end $(\mu, \nu)$ of the segment $F$ is the right vetrex of $\Delta_h$ with the minimal positive ordinate. We claim that under above assumptions, the polygon $\mathcal{N}^{(\mathbf{r})}_{\frac{\partial h}{\partial Y}}$ has at least one non-standard segment if and only if there exists $(\alpha, \beta) \in \operatorname{supp} h$ such that $0 < \beta < \nu$. Moreover:*

(a) *For every non-standard segment $T \in \mathcal{N}^{(\mathbf{r})}_{\frac{\partial h}{\partial Y}}$*

$$\frac{|T_1|}{|T_2|}\sigma(T) \leq \frac{|F_1|}{|F_2|}\sigma(F) .$$

(b) *If $T \in \mathcal{N}^{(\mathbf{r})}_{\frac{\partial h}{\partial Y}}$ is a non-standard segment parallel to $F$ then*

$$\operatorname{in}\left(\frac{\partial h}{\partial Y}, T\right) = \frac{\partial}{\partial Y}\operatorname{in}(h, F) .$$

Before the proof let see an

EXAMPLE. Let $h = X^2Y^2 + X^7 + XY^6 + X^8Y^3 + X^3Y^9 + X^9Y^6 + X^6Y^9$.

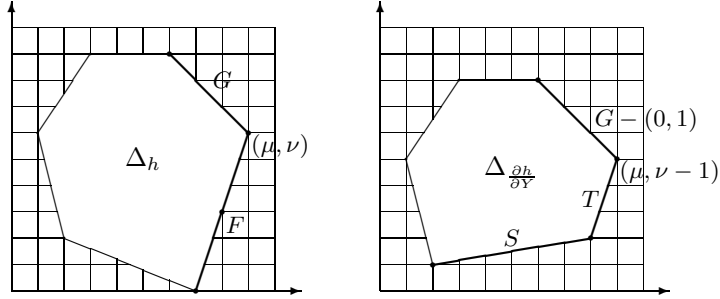

The right polygon of $h$ consists of the two segments $F$, $G$ which join the verticies $(7,0)$, $(9,6)$ and $(6,9)$, respectively. The first right vertex of $\Delta_h$ with the minimal positive ordinate is $(\mu, \nu) = (9, 6)$. There are two points of $\operatorname{supp} h$ which satisfy the condition $\{0 < \beta < \nu\}$. We have $\frac{\partial h}{\partial Y} = 2X^2Y + 6XY^5 + 3X^8Y^2 + 9X^3Y^8 + 6X^9Y^5 + 9X^6Y^8$. The righ polygon of $\frac{\partial h}{\partial Y}$ consists of the two non-standard segments $S$, $T$ and the one standard $G - (0, 1)$. Standard and non-standard segments are seperated by the vertex $(\mu, \nu - 1) = (9, 5)$. The declivities of the non standard segments are less or equal to the declivity of $F$. The segment $T$ is parallel to $F$ and obviously satisfies the condition (b). Now, let modificate the polynomial $h$ by omitting the points $(2, 2)$ and $(8, 3)$ from the support. We obtain $h = X^7 + XY^6 + X^3Y^9 + X^9Y^6 + X^6Y^9$.

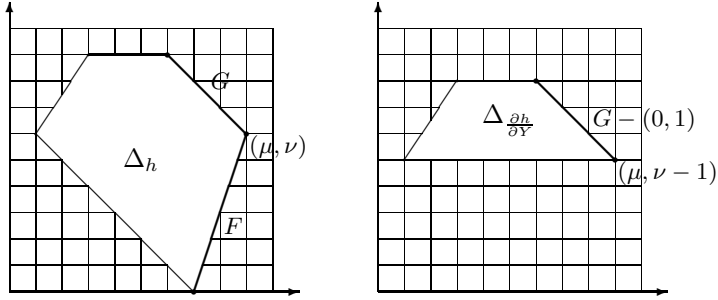

The right polygon of $h$ is exactly the same as earlier, but the right polygon of the derivative $\frac{\partial h}{\partial Y} = 6XY^5 + 9X^3Y^8 + 6X^9Y^5 + 9X^6Y^8$ has the only one standard segment $G - (0, 1)$. The reason of a lack of non-standard segments is a lack of points of $\operatorname{supp} h$ which satisfy $\{0 < \beta < \nu\}$.

Proof of Proposition 8.1. Let $T$ be a non-standard segment of $\mathcal{N}^{(\mathbf{r})}_{\frac{\partial h}{\partial Y}}$. It lies in the strip $\{0 \leq \beta \leq \nu - 1\}$, hence the lower end of the segment $T + (0, 1)$ belongs to $\operatorname{supp} h$ and satisfies $\{0 < \beta < \nu\}$. To prove the opposite implication, assume that the set $\operatorname{supp} h \cap \{0 < \beta < \nu\}$ is nonempty. Then, the set $\operatorname{supp} \frac{\partial h}{\partial Y} \cap \{\beta < \nu - 1\}$ is nonempty, too. It means that the ordinate of the right vertex $(\mu, \nu - 1)$ of $\Delta_{\frac{\partial h}{\partial Y}}$ is not minimal in $\operatorname{supp} \frac{\partial h}{\partial Y}$. Hence, there exist a segment $T \in \mathcal{N}^{(\mathbf{r})}_{\frac{\partial h}{\partial Y}}$ with the upper end $(\mu, \nu - 1)$, and clearly it is non-standard. This ends the proof of the implication. In order to prove (a) consider a non-standard segment $T$ of $\mathcal{N}^{(\mathbf{r})}_{\frac{\partial h}{\partial Y}}$. We can assume that $T$ has the maximal declivity. It means that $(\mu, \nu - 1)$ is its upper end. Let us notice that the lower end of $F$ maximizes the function $(\alpha, \beta) \mapsto \frac{\alpha - \mu}{\nu - \beta}$ in the set $\operatorname{supp} h \cap \{\beta < \nu\}$ and the maximum is equal to the declivity of $F$. Let $(\alpha_0, \beta_0)$ be the lower end of $T$. Since $(\alpha_0, \beta_0 + 1) \in \operatorname{supp} h$, we have

$$\frac{|T_1|}{|T_2|}\sigma(T) = \frac{\alpha_0 - \mu}{(\nu - 1) - \beta_0} = \frac{\alpha - \mu}{\nu - (\beta_0 + 1)} \leq \frac{|F_1|}{|F_2|}\sigma(F) \ .$$

The segment $T$ has the maximal declivity in the set of the all non-standard segments of $\mathcal{N}^{(\mathbf{r})}_h$. So, (a) is proved. In order to prove (b) let us notice that if $T \in \mathcal{N}^{(\mathbf{r})}_{\frac{\partial h}{\partial Y}}$ is a non-standard segment parallel to $F$, then it has the maximal declivity and then $(\mu, \nu - 1)$ is its upper end. In this case $T + (0, 1) \subset F$ and

$$\operatorname{in}(h, F) = \sum_{(\alpha, \beta) \in T + (0,1)} h_{\alpha\beta} X^\alpha Y^\beta \ + \ h_{p,0} X^p \ ,$$

where $(p,0)$ is the lower end of $F$. We end the proof by differentiating the formula with respect to $Y$.

A straightforward consequence of Proposition 8.1 is the following

COROLLARY 8.2. *Let* $h(X,Y) \in \mathbf{C}[X,Y]$ *and let* $T \in \mathcal{N}^{(\mathbf{r})}_{\frac{\partial h}{\partial Y}}$. *Then:*

(a) *If* $T$ *is parallel to a segment* $S \in \mathcal{N}^{(\mathbf{r})}_h$, *then*

$$\mathrm{in}\left(\frac{\partial h}{\partial Y}, T\right) = \frac{\partial}{\partial Y}\mathrm{in}(h, S),$$

(b) *If* $T$ *is not parallel to any segment from* $\mathcal{N}^{(\mathbf{r})}_h$, *then the declivity of* $T$ *is less than all the declivities of segments* $S \in \mathcal{N}^{(\mathbf{r})}_h$, *but such* $T$ *can exists only if the first right vertex of* $\Delta_h$ *lies on the horizontal axis.*

Now, we describe a structure of the right polygon of $\frac{\partial h}{\partial X}$. Notice, that $\mathrm{supp}\,\frac{\partial h}{\partial X}$ is the image of $\mathrm{supp}\,h \setminus \{\alpha = 0\}$ in the translation $(\alpha, \beta) \mapsto (\alpha - 1, \beta)$. We say that a segment $R \in \mathcal{N}^{(\mathbf{r})}_{\frac{\partial h}{\partial X}}$ is *standard* if it is the image of a segment $S \in \mathcal{N}^{(\mathbf{r})}_h$, in this translation. We write $R = S - (1,0)$. Obviously

$$\mathrm{in}\left(\frac{\partial h}{\partial X}, R\right) = \frac{\partial}{\partial X}\mathrm{in}(h, S),$$

in this case. As earlier, a necessary condition for $\mathcal{N}^{(\mathbf{r})}_{\frac{\partial h}{\partial X}}$ to be nonempty, is $\mathcal{N}^{(\mathbf{r})}_h \neq \emptyset$. For any right vertex of $\Delta_h$ with a positive abscissa, the image of this vertex in the translation $(\alpha, \beta) \mapsto (\alpha - 1, \beta)$ is a right vertex of $\Delta_{\frac{\partial h}{\partial X}}$. Similarly, if $S \in \mathcal{N}^{(\mathbf{r})}_h$ does not touch the vertical axis, then $R = S - (1,0)$ is a standard segment of $\mathcal{N}^{(\mathbf{r})}_{\frac{\partial h}{\partial X}}$. Assume that $\frac{\partial h}{\partial X} \neq 0$. Let $(\mu_1, \nu_1)$ be the first right vertex of $\Delta_h$ with a positive abscissa and let $(\mu_2, \nu_2)$ be the last one with this property. We have $\nu_1 \leq \nu_2$. Clearly, every segment $R \in \mathcal{N}^{(\mathbf{r})}_{\frac{\partial h}{\partial X}}$ in the strip $\{\nu_1 \leq \beta \leq \nu_2\}$ is standard and every segment in one of the two strips $\{0 \leq \beta \leq \nu_1\}$, $\{\nu_2 \leq \beta \leq \deg_Y h\}$, is a non-standard one. So, we can introduce, in natural way, lower and upper non-standard segments of the considered polygon. A necessary condition for existing lower non-standard segments of $\mathcal{N}^{(\mathbf{r})}_{\frac{\partial h}{\partial X}}$ is that the first right vertex of $\Delta_h$ has to lay on the vertical axis. Analogous necessary condition for existing upper non-standard segments is: the last right vertex of $\Delta_h$ has to lay on the vertical axis. In the two following propositions we give sufficient conditions for existing non-standard segments of the both types. The proofs are analogous to the proof of Proposition 8.1.

PROPOSITION 8.3. *Suppose that* $\mathcal{N}^{(\mathbf{r})}_h$ *is nonempty and the first right vertex of* $\Delta_h$ *lies on the vertical axis. Let* $F$ *be the first segment of* $\mathcal{N}^{(\mathbf{r})}_h$. *Then, the polygon* $\mathcal{N}^{(\mathbf{r})}_{\frac{\partial h}{\partial X}}$ *has at least one lower non-standard segment if and* only if there exists a point of supp $h$ with positive abscissa in the interior of the strip $\mathbf{R} \times F_2$ [1]. Moreover:

(a) For every lower non-standard segment $R \in \mathcal{N}^{(\mathbf{r})}_{\frac{\partial h}{\partial X}}$ we have

$$\frac{|R_1|}{|R_2|}\sigma(R) \leq \frac{|F_1|}{|F_2|}\sigma(F) .$$

(b) If $R \in \mathcal{N}^{(\mathbf{r})}_{\frac{\partial h}{\partial X}}$ is a lower non-standard segment parallel to $F$, then

$$\mathrm{in}\left(\frac{\partial h}{\partial X}, R\right) = \frac{\partial}{\partial X}\mathrm{in}(h, F) .$$

PROPOSITION 8.4. *Suppose that $\mathcal{N}^{(\mathbf{r})}_h$ is nonempty and the last right vertex of $\Delta_h$ lies on the vertical axis. Let $L$ be the last segment of $\mathcal{N}^{(\mathbf{r})}_h$. Then, the polygon $\mathcal{N}^{(\mathbf{r})}_{\frac{\partial h}{\partial X}}$ has at least one upper non-standard segment if and only if there exists a point of supp $h$ with positive abscissa in the interior of the strip $\mathbf{R} \times L_2$. Moreover:*

(a) *For every upper non-standard segment $R \in \mathcal{N}^{(\mathbf{r})}_{\frac{\partial h}{\partial X}}$ we have*

$$\frac{|L_1|}{|L_2|}\sigma(L) \leq \frac{|R_1|}{|R_2|}\sigma(R) .$$

(b) *If $R \in \mathcal{N}^{(\mathbf{r})}_{\frac{\partial h}{\partial X}}$ is an upper non-standard segment parallel to $L$, then*

$$\mathrm{in}\left(\frac{\partial h}{\partial X}, R\right) = \frac{\partial}{\partial X}\mathrm{in}(h, L) .$$

Let consider an

EXAMPLE. Let $h = Y + XY^2 + X^4Y^3 + Y^9 + X^3Y^7 + X^8Y^5 + X^8Y^7$.

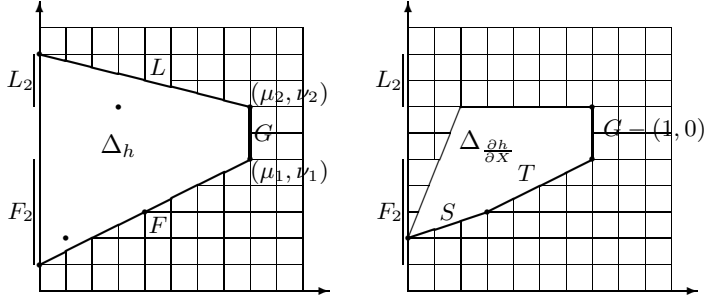

---
[1] $F_2$ is the projection of $F$ on the vertical axis

The right polygon of $\mathcal{N}_h$ consists of the three segments $F$, $G$ and $L$ which join the verticies $(0,1)$, $(8,5)$ and $(0,9)$, respectively. The right polygon of the derivative $\frac{\partial h}{\partial X} = Y^2+4X^3Y^3+3X^2Y^7+8X^7Y^5+8X^7Y^7$ consists of the the two lower non-standard segments $S$, $T$ and the one standard segment $G - (1,0)$. The reason for the lack of upper non-standard segments is the lack of points from $\operatorname{supp} h$ with positive ordinates in the interior of the strip $\mathbf{R} \times L_2$.

A straightforward consequence of Propositions 8.3 and 8.4 is the following

COROLLARY 8.5. *Let $h \in \mathbf{C}[X,Y]$ and let $R \in \mathcal{N}_{\frac{\partial h}{\partial X}}^{(\mathbf{r})}$ (notice that $\mathcal{N}_h^{(\mathbf{r})}$ is nonempty). Then:*

(a) *If $R$ is parallel to a segment $S \in \mathcal{N}_h^{(\mathbf{r})}$, then*

$$\operatorname{in}\left(\frac{\partial h}{\partial X}, R\right) = \frac{\partial}{\partial X}\operatorname{in}(h, S) \ .$$

(b) *If $R$ is not parallel to any segment of $\mathcal{N}_h^{(\mathbf{r})}$, then exactly one of the following possibilities holds:*

   (i) *The declivity of $R$ is less then all the declivities of the segments from $\mathcal{N}_h^{(\mathbf{r})}$, but such $R$ can exists only if the first right vertex of $\Delta_h$ lies on the vertical axis.*

   (ii) *The declivity of $R$ is greater then all the declivities of the segments from $\mathcal{N}_h^{(\mathbf{r})}$, but such $R$ can exists only if the last right vertex of $\Delta_h$ lies on the vertical axis.*

A consequence of the above considerations is the following theorem on a nondegeneracy inheritance.

THEOREM 8.6. *Consider a polynomial without a constant term, which is nondegenerate on each segment of its Newton polygon at infinity. Then the pair of its derivatives is also nondegenerate at infinity.*

Proof. Let $h \in \mathbf{C}[X,Y]$ be as in the statement of the theorem. Let $R$ be an arbitrary segment of the polygon of $\frac{\partial h}{\partial X}$ at infinity, and let $T$ be an arbitrary segment of the polygon of $\frac{\partial h}{\partial Y}$ at infinity. We have to show that at least one of the conditions from the definition of a pair-nondegeneracy is satisfied (Section 3). Without loss of generality, we can assume that the segments $R$ and $T$ are parallel, and both belongs to the right polygons. An existence of these segments implies that $\mathcal{N}_h^{(\mathbf{r})}$ is nonempty. We claim that there exist a segment $S \in \mathcal{N}_h^{(\mathbf{r})}$ parallel to $R$ and $T$. If not, then by Corollary 8.3 the segment $T$ has the declivity less than all the declivities of the segment from $\mathcal{N}_h^{(\mathbf{r})}$. Moreover, the first right vertex of $\Delta_h$ lies on the horizontal axis. Since $h$ is without constant term, this vertex does not lay on the vertical axis. It means, by Corollary 8.6, that the declivity of $R$ is greater then the declivities of all the segments from $\mathcal{N}_h^{(\mathbf{r})}$, what is in the contradiction with the parallelism of $R$ and $T$. So, such segment $S$ has to exist. From Corollaries 8.2 and 8.5 follows that

$$\operatorname{in}\left(\frac{\partial h}{\partial X}, R\right) = \frac{\partial}{\partial X}\operatorname{in}(h, S) \quad \text{and} \quad \operatorname{in}\left(\frac{\partial h}{\partial Y}, T\right) = \frac{\partial}{\partial Y}\operatorname{in}(h, S).$$

That means that the system $\operatorname{in}(\frac{\partial h}{\partial X}, R) = \operatorname{in}(\frac{\partial h}{\partial Y}, T) = 0$ has no solution in $(\mathbf{C} \setminus \{0\}) \times (\mathbf{C} \setminus \{0\})$ if the system $\frac{\partial}{\partial X}\operatorname{in}(h, S) = \frac{\partial}{\partial Y}\operatorname{in}(h, S) = 0$ has not. This ends the proof of the theorem.

## 9 Proof of the main result

In this section we end the proof of the main result (Theorem 2.1). Considering Theorems 3.1, 3.2 and 8.6, all we have to do is to show that the minimum of the six quantities

$$\deg \operatorname{grad} h(X, 0), \quad \inf_{R \in \mathcal{N}^{(\mathbf{r})}_{\frac{\partial h}{\partial X}}} \alpha(R, \Delta_{\frac{\partial h}{\partial Y}}), \quad \inf_{T \in \mathcal{N}^{(\mathbf{r})}_{\frac{\partial h}{\partial Y}}} \alpha(T, \Delta_{\frac{\partial h}{\partial X}}), \qquad (20)$$

$$\deg \operatorname{grad} h(0, Y), \quad \inf_{R \in \mathcal{N}^{(\mathbf{t})}_{\frac{\partial h}{\partial X}}} \beta(R, \Delta_{\frac{\partial h}{\partial Y}}), \quad \inf_{T \in \mathcal{N}^{(\mathbf{t})}_{\frac{\partial h}{\partial Y}}} \beta(T, \Delta_{\frac{\partial h}{\partial X}}), \qquad (21)$$

is equal to the minimum from the statement of Theorem 2.1, providing that a polynomial $h \in \mathbf{C}[X, Y]$ satisfies the assumptions of the theorem. Let start from a simple

PROPOSITION 9.1. *Let $h \in \mathbf{C}[X, Y]$ be a non-zero polynomial not divisible by $Y^2$. Then the right polygon of $h$ has a segment which is not exceptional if and only if $\deg_Y h \geq 2$.*

Proof. Let $S$ be a segment of the right polygon of $h$ which is not exceptional. The ordinate of the upper end of $S$ is greater than 1. So, $\deg_Y h \geq 2$. To prove the opposite implication, consider the right vertex of $\Delta_h$ with the maximal ordinate. Since $\deg_Y h \geq 2$, the ordinate of this vertex is also $\geq 2$. This vertex is the upper and of a segment, because from the assumption that $h$ is not divisible by $Y^2$ follows that there exist a point in $\operatorname{supp} h$ with the ordinate 0 or 1. Obviously, this segment is not exceptional.

By symmetry we obtain

PROPOSITION 9.2. *Let $h \in \mathbf{C}[X, Y]$ be a non-zero polynomial not divisible by $X^2$. Then the top polygon of $h$ has a segment which is not exceptional if and only if $\deg_X h \geq 2$.*

The fact that the minimum of the six quantities from (20) and (21) is equal to the minimum from the statement of Theorem 2.1, can be derived by using the two following lemmas.

LEMMA 9.3. *Let $h \in \mathbf{C}[X,Y]$ be a polynomial without constant term, such that $\frac{\partial h}{\partial X}$ is non-zero, and $h$ is not divisible by $Y^2$. Suppose that the right polygon of $h$ contains at least one segment which is not exceptional and denote by $F$ the first such segment. Let $m^{(\mathbf{r})}$ be the minimum of the three numbers from (20). Then*

$$\alpha(F) - 1 \geq m^{(\mathbf{r})} \geq \min\left\{\alpha(F) - 1,\ \alpha(F) - \frac{|F_1|}{|F_2|}\sigma(F)\right\}\ .$$

LEMMA 9.4. *Let $h \in \mathbf{C}[X,Y]$ be a polynomial without constant term, such that $\frac{\partial h}{\partial Y}$ is non-zero, and $h$ is not divisible by $X^2$. Suppose that the polygon $\mathcal{N}_h^{(\mathbf{t})}$ contains at least one segment which is not exceptional and denote by $G$ the first such segment. Let $m^{(\mathbf{t})}$ be the minimum of the three numbers from (21). Then*

$$\beta(G) - 1 \geq m^{(\mathbf{t})} \geq \min\left\{\beta(G) - 1,\ \beta(G) - \frac{|G_2|}{|G_1|}\sigma(G)\right\}$$

Obviously, Lemma 9.4 can be obtain from the previous one by symmetry. Befor a proof of Lemma 9.3 let see how the lemmas implies the desired equality.

Suppose that $h \in \mathbf{C}[X,Y]$ satisfies the assumptions of Theorem 2.1, i.e., $h$ is without constant term, with non-zero components of $\operatorname{grad} h$, not divisible by $X^2$ and not divisible by $Y^2$. Moreover, at least one of the polygons $\mathcal{N}_h^{(\mathbf{r})}$, $\mathcal{N}_h^{(\mathbf{t})}$ has a segment which is not exceptional. Without loss of generality we can assume that the right polygon has such a segment. Denote by $F$ the first segment of $\mathcal{N}_h^{(\mathbf{r})}$ which is not exceptional. We consider, further, the two possibilities: $\deg_X h \leq 1$ or $\deg_X h \geq 2$.

If $\deg_X h \leq 1$, then, by Proposition 9.2, the top polygon of $h$ is empty or has the exceptional segment, only. So, it is sufficient to show that the minimum of the six quantities from (20) and (21) is equal to $\alpha(F) - 1$. We have $|S_2| \geq 1$, for any segment $S$ in the right polygon. So, $|F_2| \geq 1$. From the assumption $\deg_X h \leq 1$ follows that $|F_1| \leq 1$. Hence, by Lemma 9.3, the minimum of the three quantities in (20) is equal to $\alpha(F) - 1$. Notice that $h = a(Y)X + b(Y)$, and since $\frac{\partial h}{\partial X} \neq 0$, we have $a \neq 0$. One can easily verify that $\deg \operatorname{grad} h(X,0) \leq 1$, what gives $\alpha(F) - 1 \leq 1$. To end the proof, in the considered case, if is enough to show that each of the three quantities in (21) is $\geq 1$. By Proposition 9.1 we have $\deg_Y h = \max\{\deg a, \deg b\} \geq 2$, hence $\deg \operatorname{grad} h(0,Y) \geq 1$. Since the top polygon of $\frac{\partial h}{\partial X}$ is empty, the second quantity is equal to $+\infty$. The third quantity is also equal to $+\infty$ if the top polygon of $\frac{\partial h}{\partial Y}$ is empty. If this polygon is nonempty, then $\deg a \geq 1$, obviously. Since $\deg \frac{\partial h}{\partial X}(0,Y) = \deg a$, we have $(0, \deg a) \in \operatorname{supp} \frac{\partial h}{\partial X}$. By the definition (3) we obtain $\beta(T, \Delta_{\frac{\partial h}{\partial X}}) \geq \deg a \geq 1$ for any $T$ such that $|T_1| \neq 0$. This ends the proof in the considered case.

Now, assume that $\deg_X h \geq 2$. By Proposition 9.2 there exists a segment of $\mathcal{N}_h^{(\mathbf{t})}$ which is not exceptional. Let $G$ by the first such segment. Denote by $m_1$ the minimum of the numbers

$$\alpha(F) - 1, \quad \beta(G) - 1,$$

and by $m_2$ the minimum of the numbers

$$\alpha(F) - \frac{|F_1|}{|F_2|}\sigma(F), \quad \beta(G) - \frac{|G_2|}{|G_1|}\sigma(G).$$

By Lemmas 9.1 and 9.2 the minimum of the six quantities in (20) and (21) less or equal $m_1$ and greater or equal $\min\{m_1, m_2\}$. To end the proof it suffices to show that $m_1 \leq m_2$. All is clear if $\frac{|F_1|}{|F_2|}\sigma(F) \leq 1$ and $\frac{|G_2|}{|G_1|}\sigma(G) \leq 1$. So, we can assume that at least one of the numbers is greater than 1. Suppose, first, that $\frac{|F_1|}{|F_2|}\sigma(F) > 1$. In particular, $\sigma(F) = 1$ and $|F_1| > |F_2|$. It is enough to show that $\beta(G) - 1 \leq m_2$. Notice, that the segment $F$ belongs to the top polygon $\mathcal{N}_h^{(\mathbf{t})}$ and its declivity $\frac{|F_2|}{|F_1|}\sigma(F)$, in this polygon, is less than 1. Since $|F_1| > |F_2| \geq 1$, the segment $F$ can not be an exceptional one in $\mathcal{N}_h^{(\mathbf{t})}$. It means that the declivity of $G$ in the top polygon is less or equal to the declivity of $F$. So, $\frac{|G_2|}{|G_1|}\sigma(G) \leq \frac{|F_2|}{|F_1|}\sigma(F) < 1$, and therefore $\beta(G) - 1 \leq \beta(G) - \frac{|G_2|}{|G_1|}\sigma(G)$. Moreover, $\beta(G) \leq \beta(F)$ and $\beta(F) - 1 \geq 0$. Hence

$$\beta(G) - 1 \leq [\beta(F) - 1]\frac{|F_1|}{|F_2|}\sigma(F) = \alpha(F) - \frac{|F_1|}{|F_2|}\sigma(F).$$

We use above a simple equality $\beta(S)\frac{|S_1|}{|S_2|}\sigma(S) = \alpha(S)$ (with $S = F$), which holds for every $S$ such that $|S_2| \neq 0$. This ends the proof of the inequality $m_1 \leq m_2$ in the considered situation. The case $\frac{|G_2|}{|G_1|}\sigma(G) > 1$ can be verified, analogously. Now we give a

*Proof of Lemma 9.3.* First, we show that for any segment $S$ of the polygon $\mathcal{N}_h^{(\mathbf{r})}$

$$\alpha(S, \Delta_{\frac{\partial h}{\partial X}}) = \alpha(S) - 1 \quad \text{and} \quad \alpha(S, \Delta_{\frac{\partial h}{\partial Y}}) = \alpha(S) - \frac{|S_1|}{|S_2|}\sigma(S). \quad (22)$$

In order to prove the first equality, consider a linear form $\phi(\alpha, \beta) = \alpha + \beta\frac{|S_1|}{|S_2|}\sigma(S)$. By the definition (2) and by the inclusion $\operatorname{supp}\frac{\partial h}{\partial X} \subset \operatorname{supp} h - (1, 0)$ we have

$$\begin{aligned}\alpha(S, \Delta_{\frac{\partial h}{\partial X}}) &= \inf\phi(\operatorname{supp}\frac{\partial h}{\partial X}) \geq \inf\phi(\operatorname{supp} h) + \phi(-1, 0) = \\ &= \alpha(S, \Delta_h) - 1.\end{aligned}$$

Using an elementary fact $\alpha(S) = \alpha(S, \Delta_h)$ for $S \in \mathcal{N}_h^{(\mathbf{r})}$, we obtain the inequality "$\geq$". To prove the opposite inequality, notice that from $\frac{\partial h}{\partial X} \neq 0$ follows that at least one of the ends of $S$ has a positive abscissa. Denote it by $(\alpha, \beta)$. Since $(\alpha - 1, \beta) \in \mathrm{supp}\, \frac{\partial h}{\partial X}$, we have

$$\alpha(S, \Delta_{\frac{\partial h}{\partial X}}) \leq (\alpha - 1) + \beta \frac{|S_1|}{|S_2|} \sigma(S);.$$

We end the proof by using a simple fact that $\alpha(S) = \alpha + \beta \frac{|S_1|}{|S_2|} \sigma(S)$ for any $(\alpha, \beta) \in S$. The second equality in (22) can be proved analogously.

Now, we give a proof of the inequalities in the statement of Lemma 9.3. Let $(\mu, \nu)$ be the first right vertex of the diagram $\Delta_h$ with the positive ordinate. We consider two cases $\nu = 1$ and $\nu \geq 2$. In both the cases we show that each quantity from (20) if greater or equal than $\alpha(F) - 1$, or $\alpha(F) - \frac{|F_1|}{|F_2|} \sigma(F)$. Moreover, we show that at least one of the considered quantities is less or equal to $\alpha(F) - 1$.

Assume that $\nu = 1$. Since $F$ is the first segment of $\mathcal{N}_h^{(\mathbf{r})}$ which is not exceptional, the vertex $(\mu, \nu) = (\mu, 1)$ is the lower end of the segment $F$. Since $(\mu, 1) \in F$, we have $\alpha(F) = \mu + \frac{|F_1|}{|F_2|} \sigma(F)$. Notice that $(\mu, 0)$ is the vertex of $\Delta_{\frac{\partial h}{\partial Y}}$, therefore $\deg \frac{\partial h}{\partial Y}(X, 0) = \mu$. Hence

$$\deg \mathrm{grad}\, h(X, 0) \geq \deg \frac{\partial h}{\partial Y}(X, 0) = \mu = \alpha(F) - \frac{|F_1|}{|F_2|} \sigma(F) \,. \qquad (23)$$

Let us note a general fact that for any non-zero $f \in \mathbf{C}[X, Y]$ and for any $S$ such that $|S_2| \neq 0$ we have $\alpha(S, \Delta_f) \geq \deg f(X, 0)$. For $f(X, 0) = 0$ it is obvious. If $f(X, 0) \neq 0$, then it follows from (2) and from the fact that $(\deg f(X, 0), 0) \in \mathrm{supp}\, f$. Using this fact to $\frac{\partial h}{\partial Y}$ we obtain

$$\inf_{R \in \mathcal{N}_{\frac{\partial h}{\partial X}}^{(\mathbf{r})}} \alpha(R, \Delta_{\frac{\partial h}{\partial Y}}) \geq \mu = \alpha(F) - \frac{|F_1|}{|F_2|} \sigma(F) \,. \qquad (24)$$

In order to estimate the third quantity in (20), notice that $F - (0, 1)$ is the first segment of the right polygon of $\frac{\partial h}{\partial Y}$. That means, it is the segment with the minimal declivity in this polygon. The value $\alpha(\,\cdot\,, \Delta_{\frac{\partial h}{\partial X}})$ does not change if we substitute the sement $F - (0, 1)$ by the parallel $F$. Hence, by using the second equality in (22) we obtain

$$\inf_{T \in \mathcal{N}_{\frac{\partial h}{\partial Y}}^{(\mathbf{r})}} \alpha(T, \Delta_{\frac{\partial h}{\partial X}}) = \alpha(F, \Delta_{\frac{\partial h}{\partial X}}) = \alpha(F) - 1 \,.$$

Joining the above equality with estimates (23) and (24) we end the proof of the lemma in the considered case.

Now, assume that $\nu \geq 2$. We claim that the first right vertex of $\Delta_h$ lies on the horizontal axis, in this case. If not, then this vertex coincide with $(\mu, \nu)$, and by $\nu \geq 2$, $h$ is divisible by $Y^2$ (contradiction). So, denote by $(p, 0)$ the first vertex of $\Delta_h$. Obviously, $p = \alpha(F)$. Since $h$ is without constant term, we have $p \geq 1$. Hence, $(p-1, 0)$ is a right vertex of $\Delta_{\frac{\partial h}{\partial X}}$ and $\deg \frac{\partial h}{\partial X}(X, 0) = p - 1 = \alpha(F) - 1$. As earlier we have

$$\deg \operatorname{grad} h(X, 0) \geq \deg \frac{\partial h}{\partial X}(X, 0) = \alpha(F) - 1$$

and

$$\inf_{T \in \mathcal{N}^{(\mathbf{r})}_{\frac{\partial h}{\partial Y}}} \alpha(T, \Delta_{\frac{\partial h}{\partial X}}) \geq \deg \frac{\partial h}{\partial X}(X, 0) = \alpha(F) - 1 \;.$$

To estimate the third quantity, notice that all the declivities of all the segments from the right polygon of $\frac{\partial h}{\partial X}$ are greater or equal to the declivity of the first segment $F$ of $\mathcal{N}^{(\mathbf{r})}_h$. It follows from Corollary 8.6, because the first right vertex $(p, 0)$ of $\Delta_h$ does not lay on the vertical axis. Hence

$$\inf_{R \in \mathcal{N}^{(\mathbf{r})}_{\frac{\partial h}{\partial X}}} \alpha(R, \Delta_{\frac{\partial h}{\partial Y}}) \geq \alpha(F, \Delta_{\frac{\partial h}{\partial Y}}) = \alpha(F) - \frac{|F_1|}{|F_2|} \sigma(F) \;.$$

To end the proof of the lemma it is sufficient to show that in any case, at least one of the quantities in (20) is less or equal to $\alpha(F) - 1$. We consider two cases connected with existence of non-standard segments of the right polygon of $\frac{\partial h}{\partial Y}$. If such a segment exist, then, due to Proposition 8.2, its declivity is less or equal to the declivity of $F$. Hence

$$\inf_{T \in \mathcal{N}^{(\mathbf{r})}_{\frac{\partial h}{\partial Y}}} \alpha(T, \Delta_{\frac{\partial h}{\partial X}}) \leq \alpha(F, \Delta_{\frac{\partial h}{\partial X}}) = \alpha(F) - 1 \;.$$

In the opposite case, by the same proposition, the set $\operatorname{supp} h \cap \{0 < \beta < \nu\}$ is empty, what, considering $\nu \geq 2$, means that $\frac{\partial h}{\partial Y}(X, 0) = 0$. Hence,

$$\deg \operatorname{grad} h(X, 0) = \max \left\{ \deg \frac{\partial h}{\partial X}(X, 0), -\infty \right\} = \alpha(F) - 1 \;.$$

This ends the prove of the lemma and the prove of the main result.

## Acknowledgments

Research is partially supported by the State Committee for Scientific Research of Poland, grant No. 2 PO3A 079 08.

Kielce University of Technology
Aleja Tysiąclecia Państwa Polskiego 7
25-314 Kielce, Poland
e-mail:lenarcik@sabat.tu.kielce.pl